\documentclass[11pt,namelimits,sumlimits]{amsart}
\usepackage{citesort}
\usepackage{amssymb,amsmath}
\usepackage[mathscr]{eucal}

\usepackage{graphicx}
\usepackage{amsmath,amsthm}
\usepackage{graphics}
\usepackage{color}
\usepackage{epsfig}
\usepackage{fullpage}
\usepackage{amssymb,amsmath}
\usepackage[mathscr]{eucal}

\newcounter{mtheorem}

\newtheorem{theorem}{Theorem}[section]
\newtheorem{lemma}[theorem]{Lemma}
\newtheorem{prop}[theorem]{Proposition}
\newtheorem{corollary}[theorem]{Corollary}

\newtheorem{definition}[theorem]{Definition}

\theoremstyle{remark}

\numberwithin{equation}{section}

\newcommand{\tb}{{\underline{\tau}}}

\newcommand{\SSS}{\mathsf{S}}

\newcommand{\group}{\widehat{\mathscr{G}}}

\newcommand{\beq}{\begin{equation}}
\newcommand{\eeq}{\end{equation}}
\newcommand{\bea}{\begin{eqnarray}}
\newcommand{\eea}{\end{eqnarray}}

\newcommand{\C}{\mathbb{C}}
\newcommand{\R}{\mathbb{R}}

\newcommand{\N}{\mathbb{N}}

\newcommand{\Sph}{\mathbb{S}}
\newcommand{\ra}{\rightarrow}

\newcommand{\Xunder}{\underline{X}}

\newcommand{\cunder}{\underline{c}\,}

\newcommand{\Stilde}{\widetilde{S}}

\newcommand{\ellunder}{{\underline{\ell}}}
\newcommand{\taa}{{\underline{t}}}
\newcommand{\ttt}{{\widetilde{t}}}
\newcommand{\azero}{{\underline{a}}}
\newcommand{\xunder}{{\underline{x}}}
\newcommand{\xover}{{\overline{x}}}
\newcommand{\xboth}{{\overline{\underline{x}}}}

\newcommand{\Lcal}{{\mathcal{L}}}
\newcommand{\Rcal}{{\mathcal{R}}}
\newcommand{\Jcal}{{\mathcal{J}}}

\newcommand{\Lchi}{{\mathcal{L}_\chi}}

\newcommand{\avg}{\operatorname{avg}}
\newcommand{\sym}{{sym}}

\newcommand{\Lh}{\ensuremath{\mathcal L_h}}

\newcommand{\xx}{\ensuremath{\mathrm{x}}}
\newcommand{\yy}{\ensuremath{\mathrm{y}}}
\newcommand{\zz}{\ensuremath{\mathrm{z}}}

\newcommand{\xxtilde}{\ensuremath{\widetilde{\mathrm{x}}}}
\newcommand{\yytilde}{\ensuremath{\widetilde{\mathrm{y}}}}
\newcommand{\domPhi}{\text{Dom}_\Phi}
\newcommand{\grouptilde}{{\mathscr{G}}}  
\newcommand{\piover}{{\pi/\sqrt{2}\,m}}  
\newcommand{\twopiover}{{\sqrt{2}\,\pi/m}}  
\newcommand{\T}{\mathbb{T}}
\newcommand{\D}{\mathbb{D}}
\newcommand{\Dhat}{\widehat{\mathbb{D}}}
\newcommand{\Xhat}{\widehat{X}}
\newcommand{\Mcat}{M_{cat}}
\newcommand{\Mtor}{M_{tor}}
\newcommand{\Mhat}{\widehat{M}}
\newcommand{\Mhatcat}{\widehat{M}_{cat}}
\newcommand{\Mhattor}{\widehat{M}_{tor}}

\newcommand{\Lchilambda}{{\mathcal{L}_{\chi,\lambda}}}

\newcommand{\ghat}{\widehat{g}}
\newcommand{\Ahat}{\widehat{A}}
\newcommand{\chihat}{\widehat{\chi}}
\newcommand{\nuhat}{\widehat{\nu}}
\newcommand{\Cunder}{{\underline{C}}}
\newcommand{\Cover}{{\overline{C}}}
\newcommand{\Fcal}{{\mathcal{F}}}
\newcommand{\XXX}{\widehat{\mathsf{X}}}
\newcommand{\YYY}{\widehat{\mathsf{Y}}}
\newcommand{\ZZZ}{\widehat{\mathsf{Z}}}

\newcommand{\xbar}{\underline{\XXX}}
\newcommand{\ybar}{\underline{\YYY}}
\newcommand{\zbar}{\underline{\ZZZ}}
\newcommand{\xbartilde}{{\underline{\mathsf{X}}}}
\newcommand{\ybartilde}{{\underline{\mathsf{Y}}}}
\newcommand{\zbartilde}{{\underline{\mathsf{Z}}}}
\newcommand{\xtilde}{{\mathsf{X}}}
\newcommand{\ytilde}{{\mathsf{Y}}}

\newcommand{\vece}{\vec{e}}
\newcommand{\vecK}{\vec{K}}
\newcommand{\veceta}{\vec{\eta}}
\newcommand{\real}{\operatorname{Re}}
\newcommand{\tr}{\operatorname{tr}}
\newcommand{\arccosh}{\operatorname{arccosh}}

\begin{document}

\title[Doubling]{Minimal surfaces in the three-sphere \\
by doubling the Clifford torus}

\author[N.~Kapouleas]{Nikolaos~Kapouleas}
\address{Department of Mathematics, Brown University, Providence,
RI 02912} \email{nicos@math.brown.edu}

\author[S.-D. Yang]{Seong-Deog~Yang}

\address{Department of Mathematics, Korea University, Anam-Dong Seongbuk-gu,
Seoul, 136-713, South Korea}
\email{sdyang@korea.ac.kr  }


\date{\today}


\keywords{Differential geometry, minimal surfaces,
partial differential equations, perturbation methods}

\begin{abstract}
We construct 
embedded closed minimal surfaces in the round three-sphere
$\mathbb S^3(1)$,
resembling two parallel copies of the Clifford torus,
joined by $m^2$ small catenoidal bridges
symmetrically arranged along a square lattice of points on the torus.
\end{abstract}

\maketitle
\section{Introduction}
\label{Sintro}
\nopagebreak

\subsection*{Historical background and the general idea}
$\phantom{ab}$
\nopagebreak

An interesting general construction for minimal surfaces is motivated by
examples of minimal surfaces which resemble two copies
of a minimal surface joined together with many catenoidal bridges.
Karcher, Pinkall, and Sterling
have constructed \cite{KPS}
minimal surfaces resembling
roughly an equatorial sphere in $\Bbb S^3(1)$ which has been ``doubled'',
and the two sheets have been connected by necks
arranged at the vertices of a Platonic solid,
with the corresponding symmetry imposed.
The examples constructed this way are finitely many,
because the Platonic solids are finitely many and the size of the neck is
determined by the neck configuration (their number and positions).
Pitts and Rubinstein have discussed \cite{PRu}
constructions for families of minimal surfaces,
where the size of the catenoidal bridges used can be arbitrarily small
and the genus then tends to infinity,
while the surfaces tend to a limit varifold.
These constructions are highly symmetrical.
Some of the constructions have a limit varifold which is a minimal
surface counted with multiplicity two.

We call such constructions ``doubling constructions'' as suggested in
\cite{kapouleas:survey}.
The ingredients for such a construction would be a minimal surface
$\Sigma$ in a Riemannian three-manifold,
two nearby copies of $\Sigma$, $\Sigma_1$ and $\Sigma_2$,
and a set of points $L\subset \Sigma$.
$\Sigma_1$ and $\Sigma_2$ can be thought of as the graphs of two functions
$\phi_1$ and $\phi_2$ on $\Sigma$.
$\phi_1$ and $\phi_2$ are assumed to be small and with small derivatives.
The minimal surface constructed would consist of a region $M_\Sigma$
which approximates $\Sigma_1$ and $\Sigma_2$ minus small discs,
and a collection of regions which approximate small truncated catenoids.
The discs removed are centered at the points on $\Sigma_1$
and $\Sigma_2$ corresponding to the points of $L$.
The catenoidal regions serve as bridges connecting to $M_\Sigma$
at the boundaries of the removed discs.
We call directions perpendicular to $\Sigma$ ``vertical'',
and directions along $\Sigma$ ``horizontal''.
The axes of the catenoidal regions would be approximately vertical.

Since a Riemannian manifold at small scale is approximately Euclidean,
we can use horizontal and vertical (approximate)
translations to find
balancing obstructions to the existence of such surfaces.
More precisely we can consider the force $F$ exerted by the region close
to $\Sigma_1$ to a catenoidal bridge,
and the force $F_c$ exerted
through the waist of the bridge
to the part of the catenoidal bridge closer to $\Sigma_1$,
by the other part.
The vertical component of $F_c$ is approximately equal to the length of its waist.
(Balancing for minimal surfaces 
is based simply on the first variation formula \cite{KKS,Si}.
For a general discussion see \cite{kapouleas:survey}.)
If $F$ is intercepted at a suitable curve which can be approximated
by a curve on $\Sigma_1$ enclosing a domain $\Omega\subset\Sigma_1$,
then the vertical component $F$ can be approximated by 
the integral of the mean curvature of $\Sigma_1$ on $\Omega$.
Because of the smallness assumptions for $\phi_1$,
we can ignore the nonlinear terms and the derivatives,
and then the mean curvature is approximated by
$(|A|^2+Ric(\nu,\nu))\phi_1$,
and the vertical component of $F$ by the area of $\Omega$ times
the value of
$(|A|^2+Ric(\nu,\nu))\phi_1$ at the corresponding point of $L$.

The above heuristic argument suggests that a necessary condition for a doubling construction
is that 
the mean curvature of the parallel surfaces to $\Sigma$ points away from $\Sigma$,
which in general amounts to
\addtocounter{theorem}{1}
\begin{equation}
\label{Econdition}
|A|^2+Ric(\nu,\nu) >0 \qquad\text{ on }\quad  \Sigma.
\end{equation}
This condition then ensures
that the vertical components of $F$
and $F_c$ point in opposite directions.
Moreover,
vertical component balancing considerations as above,
relate the size of $\phi_1$ and $\phi_2$
with the size of the catenoidal bridge and the area of $\Omega$.
Since the matching of the catenoidal bridge to $\Sigma_1$
and $\Sigma_2$ gives further relations between
$\phi_1$, $\phi_2$, and the size of the bridge,
and the area of $\Omega$ can be guessed from $L$,
it would appear that $L$ determines completely the construction.
Horizontal force considerations should further restrict the possible 
neck configurations $L$ and the sizes of the catenoidal bridges.

Developing in detail such a general construction is beyond the scope of this paper.
Instead we present a particular doubling construction where 
$\Sigma=\T$, the Clifford torus in the unit three-sphere,
and the neck configuration $L$ is a square lattice of points on $\T$.
Because of the high symmetry involved the construction simplifies significantly,
in particular we do not need to consider horizontal forces.
This construction has been outlined in \cite{kapouleas:survey}.
The method used is a gluing Partial Differential Equations method.
The particular kind of methods used relates most
closely to the methods developed in \cite{schoen,kapouleas:annals},
especially as they evolved and were systematized in \cite{kapouleas:wente}.
We refer the reader to \cite{kapouleas:survey} for a general discussion.

Another motivation for the construction in this paper is that
it is nontrivial to obtain new examples of minimal surfaces
in the round three-sphere, and the list of known
examples
is limited
\cite{L2,KPS,Ka1,Ka2}.
Moreover,
desingularization methods can be used to combine
the new surfaces produced here,
to construct a more varied class of 
further examples \cite[Theorem G]{kapouleas:survey}.
The desingularization constructions will appear elsewhere.

\subsection*{Outline of the construction}
$\phantom{ab}$
\nopagebreak

It is convenient that there is a simple coordinate system 
which is well adjusted to our purposes.
We study this coordinate system in appendix \ref{AppendixA}.
We call the corresponding coordinates $(\xx,\yy,\zz)$.
The surface $\{\zz=0\}$ in $\Sph^3(1)$
is the Clifford torus $\T$ on which the doubling construction is based.
The surfaces parallel to $\T$ are the surfaces of constant $\zz$.

In section \ref{Sinitsurf} we construct the initial surfaces $M$.
The construction is based on a square lattice $L\subset\T$ (see \ref{Elattices})
which consists of $m^2$ points.
The construction of the minimal surfaces in the main theorem
works when $m$ is large.
The surfaces constructed have genus $m^2+1$ because they amount to two tori
connected by $m^2$ handles.
The size of the catenoidal bridges $\tau$ can not be predicted precisely,
but up to a factor which is uniformly controlled independently of $m$
is given by $\tb := m^{-1}e^{-m^2/4\pi}$ (see \ref{Etau}).
This formula can be guessed from balancing considerations
as outlined above (or see \cite{kapouleas:survey}).
It allows us to prove that we can choose $\tau$ so that the construction works
(see \ref{LFcal} and the proof of the main theorem \ref{Tmain}).

The construction of $M$ is carried out in parallel with a similar construction
of a surface $\Mhat$ which would give a doubling of the plane in three-dimensional
Euclidean space.
By the maximum principle,
$\Mhat$ can not be perturbed to minimality,
in contrast with $M$ which by the main theorem of the paper \ref{Tmain}
can (for a certain $\tau$).
This is consistent with \ref{Econdition},
since $|A|^2+Ric(\nu,\nu)=4>0$
on $\T$,
while
$|A|^2+Ric(\nu,\nu)=0$
on the plane
and the mean curvature vanishes on its parallel surfaces which are planes themselves.
Actually the conormal on a perturbed $\Mhat$ on the vertical planes of reflectional
symmetry (that is on $\partial\Dhat$, see \ref{Dgroup}) is horizontal,
so the force $F$ in the discussion above vanishes
providing an alternative proof that $\Mhat$ can not be corrected.

As it is often the case in such constructions
\cite{kapouleas:annals,kapouleas:finite,kapouleas:minimal,kapouleas:imc,kapouleas:wente,kapouleas:wente:announce,kapouleas:jdg,kapouleas:cmp,kapouleas:annals,kapouleas:bulletin,kapouleas:survey,Y2,Y1},
it is convenient to define two more metrics
on the initial surfaces $M$, $h$ and $\chi$,
besides the induced metric $g$.
$h$ and $\chi$ are conformal to $g$.
$h$ allows us to write the linearized equation with uniformly bounded coefficients.
Moreover,
it allows us to understand the spectrum and the approximate kernel.
In the usual terminology $M$ modulo the symmetries has two standard regions,
which when viewed with respect to $h$ tend to a planar square and a unit sphere.
The square corresponds to a fundamental domain of $\T$
and the unit sphere to the catenoidal bridge.
There is only one (modulo the symmetries)
transition region $\Lambda$ connecting the standard regions.
$(\Lambda,\chi)$ is approximately isometric to a standard cylinder
of length $m^2/4\pi$ up to lower order terms.
The geometric quantities of $M$ are discussed in \ref{Lquantities}.
These estimates are important because they allow us to ensure that we can
perturb to minimality with an appropriately small perturbation.
Finally in \ref{Lh}
we quantify the limiting behavior of the standard regions in the 
$h$ metric as $m\to\infty$.

In section \ref{Slinear} we develop the linear theory needed.
All we need from this section is \ref{PMsolution}
and \ref{LEH}.
In \ref{LEH} we simply extract from the information we have on the mean curvature
from \ref{Lquantities} the relevant estimate we can use according
to the linear theory.
In \ref{PMsolution} we provide a solution modulo the substitute kernel
for the linear problem with appropriate decay estimates.
The construction leading to \ref{PMsolution} follows the general methodology of
\cite{kapouleas:wente}.
It is simpler than usual however,
because of the small number of standard and transition regions,
and the one-dimensionality of the substitute kernel,
which can serve also as extended substitute kernel
(see \cite{kapouleas:survey} for a general discussion).
The one-dimensionality of the approximate and (hence) the substitute kernel
follows from the fact that the symmetries kill the first harmonics
of the Laplacian on the spherical standard region corresponding to the catenoidal bridge,
and therefore the only eigenfunctions allowed in the kernel in the limiting configuration
as $m\to\infty$,
are the constants on the square (see \ref{Pappker}).
It turns out that the substitute kernel is enough for arranging
the decay we need (see \ref{Lv} and \ref{LStilde1E}),
and hence there is no need for extra ``extended substitute kernel''.

Finally in section \ref{Sresults} we prove the main theorem.
To do so we first provide in \ref{Lquadratic} an estimate of the 
nonlinear terms consistent with the decay estimates we have.
This estimate is based on a general estimate which can be derived
from general principles (see \ref{PXpert})
and which we present in appendix \ref{AppendixB}.
Next we calculate in detail the forces in the spirit of the
discussion earlier (see \ref{LFcal}),
and use that information to ensure that there is some initial surface
$M$ which 
can be perturbed to minimality.
This is consistent with the Geometric Principle (see \cite{kapouleas:wente,kapouleas:survey})
because
effectively creation of substitute kernel is achieved
by repositioning the copies of $\T$ used in the construction
at varying distances $a\tau$ from $\T$.
Finally we state and prove the main theorem \ref{Tmain}
by using as usual the Schauder fixed point theorem
\cite[Theorem 11.1]{gilbarg} to minimize the required estimates.
We remark that the minimal surfaces we find are
consistent with the description of the surfaces in Example 12
in \cite[page 306]{PRu}.

\subsection*{Notation and conventions}
$\phantom{ab}$
\nopagebreak

In this paper we use weighted H\"{o}lder norms.
A definition which works for our purposes in this paper is the following: 
\addtocounter{theorem}{1}
\begin{equation}
\label{E:weightedHolder}
\|\phi: C^{k,\beta}(\Omega,g,f)\|:=
\sup_{x\in\Omega}\frac{\,\|\phi:C^{k,\beta}(\Omega\cap B_x,g)\|\,}{f(x)},
\end{equation}
where $\Omega$ is a domain inside a Riemannian manifold $(M,g)$,
$f$ is a weight function on $\Omega$,
$B_x$ is a geodesic ball centered at $x$ and of radius the minimum of
$1$ and half the injectivity radius at $x$.

We will be using extensively cut-off functions and for this reason we adopt the
following notation:
We fix a smooth function $\Psi:\R\to[0,1]$ with the following properties:
\newline
(i).
$\Psi$ is nondecreasing.
\newline
(ii).
$\Psi\equiv1$ on $[1,\infty]$ and $\Psi\equiv0$ on $(-\infty,-1]$.
\newline
(iii).
$\Psi-\frac12$ is an odd function.
\newline
Given then $a,b\in \R$ with $a\ne b$,
we define a smooth function $\psi[a,b]:\R\to[0,1]$
by
\addtocounter{theorem}{1}
\begin{equation}
\label{E:psiab}
\psi[a,b]=\Psi\circ L_{a,b},
\end{equation}
where $L_{a,b}:\R\to\R$ is the linear function defined by the requirements $L(a)=-3$ and $L(b)=3$.

Clearly then $\psi[a,b]$ has the following properties:
\newline
(i).
$\psi[a,b]$ is weakly monotone.
\newline
(ii).
$\psi[a,b]=1$ on a neighborhood of $b$ and 
$\psi[a,b]=0$ on a neighborhood of $a$.
\newline
(iii).
$\psi[a,b]+\psi[b,a]=1$ on $\R$.

We will denote the span of vectors $e_1,\dots,e_k$ with coefficients in a field $\mathbb{F}$
by $\left< e_1,\dots,e_k\right>_\mathbb{F}$.

\subsection*{Acknowledgments}
The authors would like to thank Rick Schoen 
for his constant interest and support and insightful discussions and suggestions.
N.K. would like to thank the Mathematics Department and the MRC at Stanford University
for providing a stimulating mathematical environment and generous financial support
during Fall 2006.

\section{The initial surfaces}
\label{Sinitsurf}
\nopagebreak

In this section we define and discuss the initial surfaces.
The genus and the geometry of the initial surfaces depend on 
$m\in\N$ which we fix now and is assumed to be as large as needed.
The number of catenoidal bridges used to connect the two 
parallel copies of the Clifford torus is $m^2$
and the genus of the resulting surface $m^2+1$.
These bridges are arranged with maximal symmetry at the points
of a square lattice.
To describe the symmetry involved we have the following
(recall Appendix \ref{AppendixA}):

\addtocounter{equation}{1}
\begin{definition}
\label{Dgroup}
We denote by $\group$ the group of diffeomorphisms
of $\domPhi$ generated by the reflections
$\xbar$, $\xbar_\piover$,
$\ybar$, $\ybar_\piover$,
and
$\zbar$.
We denote by $\grouptilde$ the group of 
isometries of $\Sph^3(1)$
generated by the reflections
$\xbartilde$, $\xbartilde_\piover$,
$\ybartilde$, $\ybartilde_\piover$,
and
$\zbartilde$.
We also define $\D\subset\Sph^3(1)$ and 
$\Dhat\subset\domPhi$ by $\D:=\Phi(\Dhat)$
and
\begin{equation*}
\Dhat:=\left\{(\xx,\yy,\zz)\in\domPhi:
\,
|\xx|\le\frac{\pi}{\sqrt{2\,}\,m},\,\,
|\yy|\le\frac{\pi}{\sqrt{2\,}\,m}\right\}.
\end{equation*}
\end{definition}

The reflections
$\xbar$, $\xbar_\piover$,
$\ybar$, $\ybar_\piover$,
and
$\zbar$
generating $\group$ are with respect to the
planes $\{\xx=0\}$, $\{\yy=0\}$,
$\{\xx=\frac{\pi}{\sqrt{2\,}\,m}\}$, $\{\yy=\frac{\pi}{\sqrt{2\,}\,m}\}$,
and the line $\{\xx=\yy,\zz=0\}$ respectively.
Clearly $\XXX_\twopiover,\YYY_\twopiover\in\group$
and
$\xtilde_\twopiover,\ytilde_\twopiover\in\grouptilde$.
$\Dhat$ is a fundamental domain for the action of the translations in $\group$
and is invariant under the action of 
$\xbar$, $\ybar$ and $\zbar$.
Similarly (recall 
\ref{LPhisymmetries})
$\D$ is a fundamental domain for the action of the rotations in $\grouptilde$
and is invariant under the action of 
$\xbartilde$, $\ybartilde$ and $\zbartilde$.

We define 
square lattices $\widehat{L}$ on the plane $\{\zz=0\}$ and
$L$ on $\T$ (recall
\ref{LPhisymmetries}
and
\ref{ECT})
by
\addtocounter{theorem}{1}
\begin{equation}
\label{Elattices}
\widehat{L}:=\group(0,0,0),\qquad\qquad
{L}:=\Phi(\widehat{L})=\grouptilde\Phi(0,0,0).
\end{equation}
${L}$ consists of $m^2$ points which will be the centers of the catenoidal
bridges we use.

The size of the catenoidal bridges depends on $m$ and
on $\zeta\in\R$ which is
a parameter of the construction.
$\zeta$ is assumed to satisfy
\addtocounter{theorem}{1}
\begin{equation}
\label{Ezetarange}
|\zeta|\le \cunder,
\end{equation}
where $\cunder$ is a constant which will be chosen later.
We define then
\addtocounter{theorem}{1}
\begin{equation}
\label{Etau}
  \tb := m^{-1}e^{-m^2/4\pi},
\qquad\qquad
\tau:=e^\zeta\,\tb .
\end{equation}

We define now a constant $a>0$,
a map $\Xhat:[-a,a]\times \Sph^1\to\Dhat$,
and a truncated catenoidal bridge $\Mhatcat$ of size $\tau$,
by the following:
\addtocounter{theorem}{1}
\begin{equation}
\label{EMhatcat}
\begin{gathered}
\Mhatcat:=\Xhat([-a,a]\times \Sph^1),
\qquad
\Xhat(t,\theta):=
(r(t)\cos\theta,r(t)\sin\theta,\zz(t)),
\\
\text{ where }
\quad
r(t):=\tau\cosh t,
\quad
\zz(t):=\tau\, t,
\quad
r(a)=\frac{1}m.
\end{gathered}
\end{equation}
Note that the definition of $a$ just given,
together with \ref{Etau} and \ref{Ezetarange}
(see also \ref{Eacat}),
implies that
\addtocounter{theorem}{1}
\begin{equation}
\label{Ea}
\left|a+\zeta-\frac{m^2}{4\pi}-\log2\right|<\tb.
\end{equation}
We also define a region of a horizontal plane
(corresponding under $\Phi$ to a parallel surface to $\T$)\
together with a gluing region by
\addtocounter{theorem}{1}
\begin{equation}
\label{EMhattor}
\begin{gathered}
\Mhattor:=\left\{(\xx,\yy,\zz)\in\Dhat:\zz=\varphi(\sqrt{\xx^2+\yy^2\,}),
\,\,
\frac1m\le\sqrt{\xx^2+\yy^2}
\right\},
\\
\begin{aligned}
\text{ where }
\quad
\varphi(r):=&
\varphi_{cat}(r)
+
\psi[m^{-1},2m^{-1}](r)\,\,
(  \,      \varphi_{cat}(m^{-1})    -     \varphi_{cat}(r)  \,  ),
\\
\text{ where }
\quad
\varphi_{cat}(r):=&
\tau\arccosh \frac r \tau=
\tau\left(\log r-\log \tau+\log\left(1+\sqrt{1-\frac{\tau^2}{r^2}\,}\right)\right).
\end{aligned}
\end{gathered}
\end{equation}
Notice that $\varphi$ transits smoothly from being
$\varphi_{cat}$ in a neighborhood of $r=1/m$,
to being the constant 
\addtocounter{theorem}{1}
\begin{equation}
\label{Eacat}
\varphi_{cat}(1/m)=\tau a
\end{equation}
for $r\ge2/m$
(note that $2<\pi/\sqrt2$).
Correspondingly $\Mhattor$ extends smoothly $\Mhatcat$ close to its inner boundary circle
and transits to the plane $\zz=\varphi_{cat}(1/m)$ close to its outer boundary.
We define then smooth embedded surfaces $\Mhat\subset\domPhi$
and $\Mcat,\Mtor,M\subset\Sph^3(1)$
by 
\addtocounter{theorem}{1}
\begin{equation}
\begin{gathered}
\label{EM}
\Mhat:=\group(\Mhatcat\cup\Mhattor),
\qquad
\Mcat:=\Phi(\Mhatcat),
\qquad
\Mtor:=\Phi(\Mhattor),
\\
M:=\Phi(\Mhat)=\grouptilde(\Mcat\cup\Mtor).
\end{gathered}
\end{equation}
$\left.\Phi\right|_{\Mhat}:\Mhat\to M$ is clearly a covering map
and $M$ is a closed embedded surface of genus $m^2+1$.
We take $M$ to be our initial surface and we will prove in the Main Theorem
that for some value of $\zeta$ it can be perturbed to a nearby minimal surface.

\subsection*{Geometric quantities on the initial surfaces}
$\phantom{ab}$
\nopagebreak

We start by discussing some of the metrics we use.
We denote by $\ghat$ the standard Euclidean metric on $\domPhi$
and by $g$ the standard metric on the round sphere $\Sph^3(1)$.
Since $\Phi$ is a covering map, these metrics induce
metrics on the range and the domain of $\Phi$ respectively,
which we denote by slight abuse of notation by the same symbols.
We also use the same symbols to denote the metrics induced
on $M$, $\Mhat$ and (by using $\Xhat$) on the cylinder $\Sph^1\times[-a,a]$.
We also define cylindrical coordinates $(r,\theta,t)$
on $\Dhat$ and $\D$ by 
\addtocounter{theorem}{1}
\begin{equation}
\label{Ecoordinatesrt}
(\xx,\yy,\zz)=(r\cos\theta,r\sin\theta,\tau\,t).
\end{equation}
We define a smooth function $\rho$ on $M$ (or $\Mhat$)
by requiring it is invariant under the action of $\grouptilde$
(or $\group$) and on $\D\cap M$
(or $\Dhat\cap\Mhat$) it satisfies
\addtocounter{theorem}{1}
\begin{equation}
\label{Erho}
\rho=\frac1r
+
\psi[m^{-1},2m^{-1}](r)\,\,
\left(     \frac2m    -    \frac1r    \right).
\end{equation}
We define then smooth metrics $\chi$ and $\chihat$ on our surfaces by 
\addtocounter{theorem}{1}
\begin{equation}
\label{Echi}
\chi:=\rho^2\,g,\qquad\quad
\chihat:=\rho^2\,\ghat.
\end{equation}

We denote by $\nu$ the unit normal which satisfies $\left<\nu,\partial_\zz\right>>0$ on $\Mhattor$,
$A$ the second fundamental form induced by $g$,
by $|A|^2$ its square length,
and by $H$ the mean curvature.
We use a hat to denote the corresponding geometric quantities induced by $\ghat$.
Note that $\zz$ is constant on $\Mtor$ close to $\partial\D\cap\partial\Mtor$,
and so we will consider it extended to $M$ as a smooth function,
by requesting invariance under the action of $\grouptilde$.
We have the following:

\addtocounter{equation}{1}
\begin{lemma}
\label{Lquantities}
Assuming that $m$ is large enough in terms of $k\in\N$
the following hold:
\newline
(i).
$\|\rho^{\pm1}: C^k(M,\chihat,\rho^{\pm1})\|\le C(k)$.
\newline
(ii).
$\|\zz: C^k(M,\chihat,|\zz|+\tau)\|\le C(k)$.
\newline
(iii).
$\|\chi-\chihat: C^k(M,\chihat,|\zz|+\tau)\|\le C(k)$.
On $\Mhatcat$ we have $\chihat=dt^2+d\theta^2$.
\newline
(iv).
$\|\rho^{-2}H:C^k(M,\chi,(\tau+\rho^{-2})(|\zz|+\tau)  )\|\le C(k)$.
\newline
(v).
$\|  |A|^2-2\tau^2 \rho^4 :C^k(M,\chi,1+\tau\rho^{2})      \|\le C(k)$.
Moreover on $\Mhatcat$
we have
$\widehat{|A|^2}=2\tau^2\rho^4$. 
\end{lemma}

\begin{proof}
We first check these estimates on $\Mtor$.
$\Mtor$ is the graph of $\varphi$ (recall \ref{EMhattor})
and by the definition of $\varphi$ and \ref{Eacat} we have
\newline
(a). 
$\|\varphi-\tau a : C^k(\Mtor,  m^{2}  ( d\xx^2+d\yy^2)   \,   )   \|
\le C(k)  \,\tau$.
\newline
By \ref{Ea} we conclude
\newline
(b).
$\tau a\le m^2\tau$.
\newline
By \ref{Erho} we conclude that
\newline
(c).
$\|\rho^{\pm1} : C^k(\Mtor,  m^{2}  ( d\xx^2+d\yy^2)   \,   )   \|
\le C(k)  \,m^{\pm1}$.
\newline
By straightforward calculation we have
\begin{equation*}
\begin{aligned}
\ghat-(d\xx^2+d\yy^2)=&
\varphi^2_\xx d\xx^2+ 2\varphi_x \varphi_y d\xx d\yy + \varphi_\yy^2 d\yy^2,
\\
g-\ghat=&
\sin2\varphi \, ( d\xx^2-  d\yy^2).
\end{aligned}
\end{equation*}
(a), (b), and (c) imply then
(i), (ii), (iii), and also
\newline
(d).
$\|g-(d\xx^2+d\yy^2): C^k(\Mtor,\chi,|\zz|)\|\le C(k) $.

Using \ref{EChr} and calculating further we conclude
\begin{multline*}
\|    \,   A-(\Gamma^3_{11}d\xx^2+\Gamma^3_{22}d\yy^2) \, :
 C^k(\Mtor,  m^{2}  ( d\xx^2+d\yy^2)   \,   )   \|
\le
\\
C\,
\|(\varphi_x,\varphi_y):
 C^{k+1}(\Mtor,  m^{2}  ( d\xx^2+d\yy^2)   \,   )   \| ,
\end{multline*}
which by (a), (b) and \ref{EChr} implies that 
\newline
(e).
$\|A+d\xx^2-d\yy^2\,:C^k(\Mtor,  m^{2}  ( d\xx^2+d\yy^2)   \,   )   \|
\le C(k)\tau$.
\newline
(d) and (e) imply then (iv) and (v).

It remains to check that the estimates hold on $\Mcat$.
For convenience we adopt the notation $O(f)$ to denote
a function (or tensor field) which satisfies for each disc $D\subset\Mcat$
of radius $1$ with respect to the $\chihat$ metric the inequality
$$
\|O(f):C^k(D,\chihat)\|\le C(k)\,\|f:C^k(D,\chihat)\|.
$$
By a straightforward calculation we have then 
$$
\begin{aligned}
\Xhat_t=&\tau(\sinh t\cos \theta,   \sinh t \sin\theta, 1),
\\
\Xhat_\theta=&
\tau\cosh t(-\sin\theta,   \cos \theta, 1),
\end{aligned}
$$
which implies
\newline
(f).
$\ghat=r^2(dt^2+d\theta^2)$
and
$g=r^2(dt^2+d\theta^2+O(\zz)\,  )$.
\newline
Using (f) and the definitions, (i), (ii), and (iii) follow.
Using \ref{EChr} and calculating further we find that
\newline
(g).
$\Ahat=\tau(-dt^2+d\theta^2)$,
and
$A=\Ahat+\widetilde{A}+O(\,(\tau+r^2)\zz)$,
where
\newline
$\widetilde{A}=r^2\,(\cos2\theta\,(-dt^2+d\theta^2)+2\sin2\theta \,dt\,d\theta\,)$.
\newline
Using (f) and (g) it is straightforward to check (iv) and (v) and complete the proof
(notice also that $\tau^2\zz< C\,\tau\, r^2$).
\end{proof}

\subsection*{Standard and transition regions}
$\phantom{ab}$
\nopagebreak

We proceed to define carefully the 
various regions on the initial surface $M$ in the usual fashion of 
\cite{kapouleas:annals,kapouleas:finite,kapouleas:minimal,kapouleas:imc,kapouleas:wente,kapouleas:wente:announce,kapouleas:jdg,kapouleas:cmp,kapouleas:annals,kapouleas:bulletin,kapouleas:survey}.
Modulo the symmetries imposed,
there are only two standard regions which we denote
by $S[0]$ (corresponding to the catenoidal bridge)
and $S[1]$ (corresponding to the torus),
and only one transition region we denote by $\Lambda$.
The extended standard regions
$\Stilde[0]$ and $\Stilde[1]$
are the standard regions augmented by the transition region.

In order to ensure uniformity with respect to different values 
of the parameter $\zeta$,
we define and use a variant $\taa$ of the parameter
$t$ by
\addtocounter{theorem}{1}
\begin{equation}
\label{Etaa}
\taa=\frac{\,\azero\,}{a}\,t,
\end{equation}
where $a$ is defined as in \ref{EMhatcat}
for the current value of $\zeta$
and $\azero$ is defined in the same way when $\zeta=0$ and hence $\tau=\tb$ (recall \ref{Etau}).
This way the range of values of $\taa$ on $\Mcat$
is $[-\azero,\azero]$ and it does not depend on $\zeta$ or $\tau$.
Note also that by \ref{Ea} we have
\addtocounter{theorem}{1}
\begin{equation}
\label{Ettaa}
|t-\taa|\le C\cunder \text{ on } \Stilde[0],
\qquad\qquad
\left|\frac {\,\azero\,}{a}-1\right|\le C \cunder m^{-2}.
\end{equation}

We use a constant $b$ to control the exact size of the standard and transition regions.
$b$ will be determined later so that the linearized equation and its spectrum
behave appropriately.
We use subscripts $x$ and $y$ to modify the usual sizes and boundary circles.
In particular each $S_x[n]$ is a neighborhood of $S[n]$,
while $\Stilde_x[n]$ is $\Stilde[n]$
with an appropriate neighborhood of its boundary excised.

\addtocounter{equation}{1}
\begin{definition}
\label{D:regions}
We define the following:
\addtocounter{theorem}{1}
\begin{subequations}
\label{E:regions}
\begin{align} 
S_x[0] &:= M \cap \D \cap \{\taa\in[-b-x,b+x]\}, \\
S_x[1] &:= M \cap \D \cap \{\taa\ge \azero-b-x\}, \\
\Stilde_x[0] &:= M \cap \D \cap \{\taa\in[-\azero+b+x,\azero-b-x]\}, \\
\Stilde_x[1] &:= M \cap \D \cap \{\taa\ge b+x\}, \\
\label{ELambda}
\Lambda_{x,y} &:= M \cap \D \cap \{\taa\in[ b+x,a-b-y]\}, \\
C_x[0] &:= M \cap \D \cap \{\taa=b+x]\}, \\
C_x[1] &:= M \cap \D \cap \{\taa= \azero-b-x\},\\
\label{ECpartial}
C_\partial &:= \partial \D \cap \partial \Mtor
\end{align}
\end{subequations}
where $b>5$ is a constant chosen finally
in the proof of \ref{LStilde1E} independently of $m$,
and 
$0\le x,y< \frac13 a-b$.
When $x=y=0$ we drop the subscripts.
We also write $\Lambda_x$ for $\Lambda_{x,x}$.
\end{definition}

The limiting behavior of the standard regions,
and the linearized operator on them,
as $m\to\infty$,
is best understood in the $h$ metric which is defined on our surfaces by
\addtocounter{theorem}{1}
\begin{equation}
\label{Eh}
h:=\frac{|A|^2+m^2}{2}g.
\end{equation}
We define the map $\varpi:\D\to\R^2$ by
\addtocounter{theorem}{1}
\begin{equation}
\label{Evarpi}
\varpi(\xx,\yy,\zz):=\frac{m}{\sqrt2}(\xx,\yy).
\end{equation}
The following lemma describes the limiting behavior as $m\to\infty$:

\addtocounter{equation}{1}
\begin{lemma}
\label{Lh}
If $m$ is large enough in terms of $b+x$, then the following hold:
\newline
(i).
$\|h-\nuhat^*g:C^5(S_x[0],\nuhat^*g)\|
\le
C(b+x)\,\tau$,
where $\nuhat^*g$ is the pullback of the standard metric of the unit sphere $\Sph^2(1)$
by $\nuhat$ and satisfies
$\nuhat^*g=\frac12 \widehat{|A|^2} \ghat=\tau^2r^{-4}\ghat= \tau^2 r^{-2}\chihat$.
Moreover $\nuhat(S_x[0])=\{(\xx,\yy,\zz)\in\Sph^2(1):\xx^2+\yy^2\ge\check{R}_x^2\}$,
where $\check{R}_x=1/\cosh[(b+x)a/\azero]$.
\newline
(ii).
$\|h-\varpi^*g:C^5(S_x[0],\varpi^*g)\|
\le
C(b+x)/m^2$,
where $\varpi^*g$ is the pullback of the standard Euclidean metric on $\R^2$
by $\varpi$ (restricted to $S_x[1]$).
Moreover there is $\widetilde{R}_x$ such that
\newline
$|\widetilde{R}_x-2^{-1/2}e^{-(b+x)a/\azero}|\le \tau$ and
$\varpi(S_x[1])=\{(\xxtilde,\yytilde)\in\R^2:
|\xxtilde|\le\frac\pi2,\, |\yytilde|\le\frac\pi2,\, \xxtilde^2+\yytilde^2\ge \widetilde{R}_x^2\}$.
\end{lemma}

\begin{proof}
Since the catenoid is a minimal surface it follows from standard theory that 
$\nuhat^*g=\frac12 \widehat{|A|^2} \ghat$,
and the expressions in terms of $r$ follow from \ref{Lquantities}.v and the definitions.
This implies that the length of $\nuhat(C_x[0])$ is $2\pi\tau\,/r(t)=2\pi/\cosh[(b+x)a/\azero]$,
which implies that $\nuhat(S_x[0])$ is as stated.
Since
$$
h-\nuhat^*g =
\frac12 (|A|^2+m^2-2\tau^2\rho^4)\,\rho^{-2}\chi+ \tau^2\rho^2(\chi-\chihat),
$$
we conclude by using \ref{Lquantities}
that
$$
\|h-\nuhat^*g:C^5(S_x[0],\chihat)\|
\le
C\,(m^2r^2+\tau)
\le
C\,\tau.
$$
This implies the desired estimate and completes the proof of (i).

The second part of (ii) follows easily from the definitions and the observation
that $\widetilde{R}_x=(m/\sqrt2)\,r(a-(b+x)a/\azero)$.
By writing
$$
h-\varpi^*g=
\frac{|A|^2+m^2}{2}(g-(d\xx^2+d\yy^2)\,)
+
\frac{|A|^2}{m^2} \varpi^*g,
$$
using \ref{Lquantities}.i to establish the analogue of (c) in the proof of \ref{Lquantities},
estimating $g-(d\xx^2+d\yy^2)$ as for (d) in the proof of \ref{Lquantities},
and estimating $|A|^2$ by \ref{Lquantities}.v,
we conclude the proof.
\end{proof}

\section{The Linearized Equation}
\label{Slinear}
\nopagebreak

\subsection*{Introduction}
$\phantom{ab}$
\nopagebreak

In this section we study the linearized equation on $M$
which can be stated in any of the following
equivalent formulations,
\addtocounter{theorem}{1}
\begin{equation}
\label{Elinear}
\Lchi u=E,
\quad \text{ or } \quad
\Lcal u=\rho^2 E,
\quad \text{ or } \quad
\Lh u=\frac{2\rho^2}{|A|^2+m^2} E,
\qquad
\end{equation}
where the corresponding linear operators are given by
\addtocounter{theorem}{1}
\begin{equation}
\label{Eoperators}
\begin{gathered}
\Lchi:=\Delta_\chi+\rho^{-2}(|A|^2+2),
\qquad
\Lh:=\Delta_h+2\frac{|A|^2+2}{|A|^2+m^2},
\\
\Lcal:=\Delta_g+|A|^2+2=\rho^2\Lchi=\frac{|A|^2+m^2}2\Lh.
\end{gathered}
\end{equation}

\subsection*{The linearized equation on the transition region}
$\phantom{ab}$
\nopagebreak

In this subsection we consider the linearized equation on the transition region
$\Lambda_{x,y}$
defined as in \ref{ELambda}, where we assume that $x,y\in[0,4]$.
For simplicity in this subsection we will denote the neck under consideration by $\Lambda$, 
and its boundary circles $C_x[0]$ and $C_y[1]$ by $\Cunder$ and $\Cover$ respectively.
We next define $\xunder,\xover,\xboth:\Lambda\ra\R$
to measure the $\taa$-coordinate distance from $\Cunder$, $\Cover$, and $\partial\Lambda=\Cunder\cup\Cover$
respectively:
\addtocounter{theorem}{1}
\begin{equation}
\label{Exdistance}
b+x+\xunder
=\taa,
\qquad
\azero-b-y-\xover
=\taa,
\qquad
\xboth:=\min(\xunder,\xover).
\end{equation}
Note that we can use $\Phi\circ\Xhat$ to identify $\Lambda$
with the cylinder $[(b+x)a/\azero,a-(b+y)a/\azero]\times\Sph^1$.
We define
$\ellunder$ to be the $\taa$-coordinate length of the cylinder
and
$\ell$ to be the $t$-coordinate length of the cylinder,
so that
\addtocounter{theorem}{1}
\begin{equation}
\label{Eell}
\ellunder=\azero-2b-x-y,
\qquad
\ell=a-(2b+x+y)a/\azero.
\end{equation}
Using \ref{Ea}
and our assumption that $x,y\in[0,4]$,
we estimate
\addtocounter{theorem}{1}
\begin{equation}
\label{Eellm}
\left|\ell+2b+\zeta-\frac{m^2}{4\pi}\right|<10.
\end{equation}

Our understanding of the linear equations on the transition region
is based on the comparison with $\Delta_\chi$,
which is based on the following lemma:

\addtocounter{equation}{1}
\begin{lemma}
\label{LLambdasmall}
The following hold on $\Lambda$:
\newline
(i).
$\|\chi-\chihat: C^5(M,\chihat)\|\le Cm^2\tau$.
\newline
(ii).
$\|\rho^{-2}(|A|^2,m^2):C^5(\Lambda,\chi,e^{-3\xboth/2})\|\le C\,e^{-3b/2}$.
\end{lemma}

\begin{proof}
This is a straightforward consequence of \ref{Lquantities}, \ref{Ea},
and the various definitions.
\end{proof}

\addtocounter{equation}{1}
\begin{prop}
\label{PLambdalow}
If $m$ is large enough then 
the lowest eigenvalue of the Dirichlet problem for $\Lchi$ on $\Lambda$ is $>C\ell^{-2}$.
\end{prop}

\begin{proof}
The proof is similar to the arguments leading to Proposition 2.28 in \cite{kapouleas:wente}.
It is easy to prove that for $\phi\in L^2(\Lambda)$ with $L^2$ derivatives and $\phi=0$ on
$\partial\Lambda$ we have
$$
\int_\Lambda 
e^{-3\xboth/2} \phi^2 d\chihat \le C\int_\Lambda |\nabla \phi|_{\chihat}^2 d\chihat,
$$
which together with 
\ref{LLambdasmall}
implies
$$
\int_\Lambda|\nabla \phi|_{\chi}^2 d\chi
-\int_\Lambda\rho^{-2}(|A|^2+2)\phi^2d\chi
\ge
(\frac23-C e^{-3b/2} )
\int_\Lambda |\nabla \phi|_{\chihat}^2 d\chihat.
$$
Using the variational characterization of eigenvalues and assuming $b$ large enough
the result follows since the smallest eigenvalue for $\Delta_{\chihat}$
is $>C\ell^{-2}$.
\end{proof}

\addtocounter{equation}{1}
\begin{corollary}
\label{CLunique}
(i).
The Dirichlet problem for $\Lchi$ on $\Lambda$ for given $C^{2,\beta}$ Dirichlet data has a unique solution.
\newline
(ii). 
For $E\in C^{0,\beta}(\Lambda)$ there is a unique $\varphi\in C^{2,\beta}(\Lambda)$
such that $\Lchi \varphi=E$ on $\Lambda$ and $\varphi=0$ on $\partial \Lambda$.
Moreover
$
\|\varphi :C^{2,\beta}(\Lambda,\chi)\|
\le
C(\beta)\,\ell^2\,
\|E:C^{0,\beta}(\Lambda,\chi)\|.
$
\end{corollary}

\begin{proof}
(i) follows trivially and (ii) by using standard linear theory.
\end{proof}

All our constructions have to respect the symmetries imposed,
in particular we only consider functions on $M$ which are invariant
under the action of $\grouptilde$.
$\Lambda$ is not invariant under $\grouptilde$ but it is invariant
under $\xbartilde$ and $\ybartilde$.
Under the identification of $\Lambda$ with a cylinder as discussed above,
$\xbartilde$ corresponds to $\theta\to\pi-\theta$,
and $\ybartilde$ corresponds to $\theta\to-\theta$.
We use the subscript ``$\SSS$'' to specify subspaces of functions on $\Lambda$
which are invariant under these symmetries.
In the next Proposition and its Corollary,
we study the Dirichlet problem
when we are allowed to modify the lowest harmonic on the boundary data
in order to have decay estimates appropriate for our purposes:

\addtocounter{equation}{1}
\begin{prop}
\label{PLambdaEV}
Assuming $b$ large enough in terms of given $\beta,\gamma\in(0,1)$,
there is a linear map $\Rcal_\Lambda:C^{0,\beta}_\SSS(\Lambda)\to C^{2,\beta}_\SSS(\Lambda)$
such that the following hold for
$E\in C^{0,\beta}_\SSS(\Lambda)$ and $V:=\Rcal_\Lambda\, E$:
\newline
(i).
$\Lchi V= E$ on $\Lambda$.
\newline
(ii).
$V$ is constant on $\Cover$ and vanishes on $\Cunder$.
\newline
(iii).
$\|V:C^{2,\beta}_\SSS(\Lambda,\chi,e^{-\gamma\xover})\|
\le
C(\beta,\gamma)\,
\|E:C^{0,\beta}_\SSS(\Lambda,\chi,e^{-\gamma\xover})\|$.
\newline
(iv).
$\Rcal_\Lambda$ depends continuously on $\tau$.

The proposition still holds if the roles of $\Cunder$ and $\Cover$ are exchanged in (ii)
and $\xover$ is replaced by $\xunder$ in (iii).
Another possibility is to allow 
$V$ to be constant on each of $\Cover$ and $\Cunder$
in (ii),
while $\xover$ is replaced by $\xboth$ in (iii).
\end{prop}

\begin{proof}
The proposition follows by standard theory if $\Lchi$ is replaced by $\Delta_{\chihat}$.
We denote the corresponding linear map and solution in the $\Delta_\chi$ case
by $\widetilde{\Rcal}_\Lambda$ and $\widetilde{V}$ respectively.
Using then \ref{LLambdasmall} we have
$$
\|\Lchi\,\widetilde{V}:C^{0,\beta}(\Lambda,\chihat,e^{-\gamma\xover})\|
\le C(\beta,\gamma)\,
(m^2\tau+ e^{-3b/2} )
\|E:C^{0,\beta}(\Lambda,\chihat,e^{-\gamma\xover})\|,
$$
and the proposition then follows by an iteration where we treat $\Lchi$ and $\Rcal_\Lambda$
as small perturbations of $\Delta_{\chihat}$ and $\widetilde{\Rcal}_\Lambda$ and assuming
$b$ and $m$ large enough.
\end{proof}

We will only need the next statement with $\varepsilon_1=1$:

\addtocounter{equation}{1}
\begin{corollary}
\label{CLambdauV}
Assuming $b$ large enough in terms of given $\beta,\gamma\in(0,1)$
and
$\varepsilon_1>0$, 
there is a linear map 
$$
\Rcal_\partial:
\{u\in C_\SSS^{2,\beta}(\Cover):\int_{\Cover}u d\theta=0\}
\to C_\SSS^{2,\beta}(\Lambda)
$$
such that the following hold for $u$ in the domain of $\Rcal_\partial$ and $V:=\Rcal_\partial u$:
\newline
(i).
$\Lchi V= 0$ on $\Lambda$.
\newline
(ii).
$V-u$ is constant on $\Cover$ and $V$ vanishes on $\Cunder$.
\newline
(iii).
$|V-u|
\le
\varepsilon_1
\,\|u:C_\SSS^{2,\beta}(\Cover,d\theta^2)\|$.
\newline
(iv).
$\|V:C_\SSS^{2,\beta}(\Lambda,\chi,e^{-\gamma\xover})\|
\le
C(\beta,\gamma)
\,\|u:C_\SSS^{2,\beta}(\Cover,d\theta^2)\|$.
\newline
(v).
$\Rcal_\partial$ depends continuously on $\tau$.

The Proposition still holds if the roles of $\Cover$ and $\Cunder$ are exchanged 
and $\xover$ is replaced by $\xunder$.
\end{corollary}

\begin{proof}
By standard theory there is a linear map
$$
\widetilde{\Rcal}_\partial:
\{u\in C_\SSS^{2,\beta}(\Cover):\int_{\Cover}u d\theta=0\}
\to C_\SSS^{2,\beta}(\Lambda)
$$
such that for $u$ in the domain and $\widetilde{V}=\widetilde{\Rcal}_\partial u$
the following hold:
\newline
(a).
$\Delta_{\chihat} \widetilde{V}= 0$ on $\Lambda$.
\newline
(b).
$\widetilde{V}=u$ on $\Cover$
and $\widetilde{V}$ vanishes on $\Cunder$.
\newline
(c).
$\|\widetilde{V}:C_\SSS^{2,\beta}(\Lambda,\chi,e^{-\gamma\xover})\|
\le
C(\beta,\gamma)
\, \|u:C_\SSS^{2,\beta}(\Cover,d\theta^2)\|$.

The corollary then follows by defining
$$
\Rcal_\partial u:=
\widetilde{\Rcal}_\partial u
-
\Rcal_\Lambda\,
\Lchi
\widetilde{\Rcal}_\partial u,
$$
applying the Proposition,
and
using
\ref{LLambdasmall}.
\end{proof}

\addtocounter{equation}{1}
\begin{corollary}
\label{CLambdau}
If $u\in C_\SSS^{2,\beta}(\Lambda)$ satisfies
$\Lchi u=0$ on $\Lambda$,
then
$$
\|u:C_\SSS^{2,\beta}(\Lambda,\chi)\|
\le
C(\beta)\,
\|u:C_\SSS^{2,\beta}(\partial\Lambda,\chi)\|.
$$
\end{corollary}

\begin{proof}
Because of 
\ref{CLambdauV} and 
\ref{CLunique}
it is enough to prove the Corollary when $u$
is constant on each boundary circle.
Let $\widetilde{V}$ be the solution of
$$
\Delta_{\chihat} \widetilde{V}= 0 \text{ on }\Lambda,
\qquad
\widetilde{V}=1 \text{ on } \Cover,
\qquad
\widetilde{V}=0 \text{ on } \Cunder.
$$
By \ref{LLambdasmall} we can write $\Lchi\widetilde{V}=E_1+E_2$
where
$\|E_1:C_\SSS^{2,\beta}(\Lambda,\chi,e^{-\gamma\xover})\|\le C\,e^{-3b/2}$
and
$\|E_2:C_\SSS^{2,\beta}(\Lambda,\chi,e^{-\gamma\xunder})\|\le C\,e^{-3b/2}/\ell$.
By applying twice 
\ref{PLambdaEV}
and assuming $b$ large enough
we obtain $\overline{V}\in C^{2,\beta}_\SSS$
such that 
$\Lchi \overline{V}= 0 \text{ on }\Lambda,$
$\overline{V}$ is constant on each boundary circle of $\Lambda$,
$\|\overline{V}:C_\SSS^{2,\beta}(\Lambda,\chi)\|
\le
C(\beta)$,
$|\overline{V}-1|\le1/9$ on $\Cover$,
and
$|\overline{V}|\le1/9$ on $\Cunder$.
By exchanging $\Cover$ with $\Cunder$ we obtain $\underline{V}$ instead of
$\overline{V}$.
By considering linear combinations of $\overline{V}$ and $\underline{V}$
we complete the proof.
\end{proof}

\subsection*{The approximate kernel}
$\phantom{ab}$
\nopagebreak

We proceed now to discuss the approximate kernel of $\Lh$ on the extended 
standard regions, \textit{cf.} \cite[Prop. 2.22]{kapouleas:wente}.
By approximate kernel we mean the span of eigenfunctions whose eigenvalues are close to $0$.
Since we have to take into account the symmetries imposed,
note that the
stabilizer of $\Stilde[0]$ with respect to the action of $\grouptilde$
is generated by the reflections 
$\xbartilde$, $\ybartilde$, and $\zbartilde$,
and the stabilizer of $\Stilde[1]$ by $\xbartilde$ and $\ybartilde$.
Therefore we have to restrict our attention to functions on the extended standard
regions which are invariant under the action of these subgroups.
Moreover the functions on $\Stilde[1]$ should extend smoothly to
$\grouptilde \Stilde[1]$.
\addtocounter{equation}{1}
\begin{definition}
\label{Dsym}
We call functions which satisfy the above conditions appropriately symmetric
and we use the subscript ``$\sym$'' to denote subspaces of appropriately symmetric functions.
\end{definition}

We understand the approximate kernel in the next proposition by comparing it to
the kernel of the operator $\Delta+2$ on the round sphere $\Sph^2(1)$,
and $\Delta$ on the square
$[-\frac\pi2,\frac\pi2]\times [-\frac\pi2,\frac\pi2]$
with Neumann boundary conditions on the boundary.
Because of the symmetries the former is trivial and the latter one-dimensional:

\addtocounter{equation}{1}
\begin{prop}
\label{Pappker}
Assuming $b$ large enough in absolute terms,
and $\tb$ small enough (equivalently $m$ large enough)
in terms of a given $\varepsilon > 0$,
the following hold:
\newline
(i).
$\Lh$ acting
on appropriately symmetric functions on $\Stilde[0]$
with vanishing Dirichlet conditions,
has no eigenvalues
in $[-1,1]$ and the corresponding approximate kernel is trivial.
\newline
(ii).
$\Lh$ acting
on appropriately symmetric
functions
on $\Stilde[1]$
has exactly one eigenvalue $\lambda_0$
in $[-\varepsilon,\varepsilon]$,
and no other eigenvalues in $[-1/2,1/2]$,
and therefore the corresponding approximate kernel is one-dimensional.
Moreover the approximate kernel is spanned by a function
$f_0\in C^\infty_{\sym}(\Stilde[1])$
which depends continuously on $\zeta$ and satisfies
$$
\Vert  f_0- 1 \,: C^{2,\beta}(S_5[1]) \Vert  < \varepsilon,
\qquad
\Vert  f_0 \,: C^{2,\beta}(\Stilde[1],\chi) \Vert  < C.
$$
\end{prop}

\begin{proof}
The proof is based on the results of \cite[Appendix B]{kapouleas:annals}
which are based on basic facts about eigenvalues and eigenfunctions \cite{chavel}.
Before using those results we remark the following:
First, the first inequality in \cite[B.1.6]{kapouleas:annals}
should read
$$
\|F_i f\|_\infty \le2\|f\|_\infty
$$
instead.
Second, the spaces of functions can be constrained to satisfy appropriate symmetries,
as indeed was the case in some of the constructions in \cite{kapouleas:annals},
and will be the case here.
Third, the only use of the Sobolev inequality
\cite[B.1.5]{kapouleas:annals}
is to establish supremum bounds for the eigenfunctions.
These in our case can be alternatively established by using
the uniformity of geometry of $S_5[n]$
to obtain interior estimates on $S_1[n]$,
and then using a variant of
\ref{CLambdau}
to obtain estimates on the transition regions.
More precisely the eigenvalue equation under consideration
is $\Lh u + \lambda u=0$, which is equivalent to
\addtocounter{theorem}{1}
\begin{equation}
\label{ELchilambda}
\Lchilambda u=0
\quad \text{ where }\quad
\Lchilambda=\Delta_\chi+
\frac{|A|^2+m^2}{2\rho^2} \lambda.
\end{equation}
Since the modified part of the operator,
$\frac{|A|^2+m^2}{2\rho^2} \lambda$,
satisfies the same estimates by \ref{LLambdasmall} (assume $|\lambda|<9$)
as $\rho^{-2}(|A|^2+2)$,
we can repeat the arguments leading to
\ref{CLambdau}
to establish the same estimate under the modified assumption that
$\Lchilambda u=0$ on $\Lambda$.

For (i) we compare with the following:
$$
N[0]=\Sph^2(1)\bigcup  \left(  \{1,-1\}\times \widetilde{D}(\widetilde{R}_0)  \right)
\quad\text{ where }\quad
\widetilde{D}(\widetilde{R}_0)=\{(\xxtilde,\yytilde)\in\R^2:\xxtilde^2+\yytilde^2\le \widetilde{R}_0^2\},
$$
where $\widetilde{R}_x$ was defined in \ref{Lh}.
The action of
$\xbartilde$, $\ybartilde$, and $\zbartilde$, on 
$N[0]$ should be consistent with their action on $M$
(recall \ref{Edomainsymmetries})
and the maps $\nuhat$ and $\varpi$:
We define for
$(\xx,\yy,\zz)\in \Sph^2(1)$ and
$(i,\xxtilde,\yytilde)\in\{1,-1\}\times\widetilde{D}(\widetilde{R}_0)$
\addtocounter{theorem}{1}
\begin{equation}
\label{Ehsymmetries}
\begin{aligned}
\xbartilde(\xx,\yy,\zz)=&(-\xx,\yy,\zz),
\qquad   &
\xbartilde(i,\xxtilde,\yytilde)=&(i,-\xxtilde,\yytilde),
\\
\ybartilde(\xx,\yy,\zz)=&(\xx,-\yy,\zz),
\qquad   &
\ybartilde(i,\xxtilde,\yytilde)=&(i,\xxtilde,-\yytilde),
\\
\zbartilde(\xx,\yy,\zz)=&(\yy,\xx,-\zz),
\qquad   &
\zbartilde(i,\xxtilde,\yytilde)=&(-i,\yytilde,\xxtilde).
\end{aligned}
\end{equation}
We consider the Dirichlet problem on $N[0]$ where the operator
is $\Delta+2$
on $\Sph^2(1)$,
and the standard Laplacian $\Delta$ on
$\{1,-1\}\times \widetilde{D}(\widetilde{R}_0)$.
By standard theory then there are no eigenvalues in $[-1,1]$ because
the symmetries do not allow the first harmonics on $\Sph^2(1)$,
and $\widetilde{R}_0$ is small enough so that the smallest eigenvalue on the discs is $>2$.

For (ii) we compare with the following:
$$
\begin{gathered}
N[1]=\check{D}\bigcup 
\left(  [-\pi/2,\pi/2]\times[-\pi/2,\pi/2]  \right),
\\
\quad\text{ where }\quad
\check{D}=\{(\xx,\yy,\zz)\in\Sph^2(1):\xx^2+\yy^2\le\check{R}_0^2, \quad \zz\ge0\},
\end{gathered}
$$
where $\check{R}_0=1/\cosh(ab/\azero)$
(recall \ref{Lh}).
The action of
$\xbartilde$ and $\ybartilde$ on
$N[1]$ should be consistent again with their action on $M$
(recall \ref{Edomainsymmetries})
and the maps $\nuhat$ and $\varpi$:
We define for
$(\xx,\yy,\zz)\in \check{D}$ and
$(\xxtilde,\yytilde)\in
[-\pi/2,\pi/2]\times[-\pi/2,\pi/2]$
\addtocounter{theorem}{1}
\begin{equation}
\label{Ehsymmetries1}
\begin{aligned}
\xbartilde(\xx,\yy,\zz)=&(-\xx,\yy,\zz),
\qquad   &
\xbartilde(\xxtilde,\yytilde)=&(-\xxtilde,\yytilde),
\\
\ybartilde(\xx,\yy,\zz)=&(\xx,-\yy,\zz),
\qquad   &
\ybartilde(\xxtilde,\yytilde)=&(\xxtilde,-\yytilde),
\end{aligned}
\end{equation}
As before the operator on $\check{D}\subset\Sph^2(1)$
is $\Delta+2$
and on 
$[-\pi/2,\pi/2]\times[-\pi/2,\pi/2]\subset\R^2$
is the standard Laplacian $\Delta$.
The boundary conditions are the Dirichlet condition on $\partial\check{D}$
and the Neumann condition---more precisely extendibility to $\R^2$
by reflections across the lines $\{\xxtilde=n\pi/2\}$
and $\{\xxtilde=n\pi/2\}$ ($n\in\N$)---for the boundary
of the square $[-\pi/2,\pi/2]\times[-\pi/2,\pi/2]\subset\R^2$.
The smallness of $\check{R}_0$
and our knowledge of the eigenvalues on the square imply
the only eigenvalue in $[-2/3,2/3]$ is $0$,
with corresponding eigenfunctions the functions which are constant on the square
and vanish on $\check{D}$.

To complete the proof we use $\nuhat$, $\varpi$, and the logarithmic cut-off 
function $\psi[2d,d]\circ\xboth$ on $\Lambda$
to define the maps $F_1$ and $F_2$
required by \cite[B.1.4]{kapouleas:annals}
as usual.
$d$ is taken to be large enough in terms of 
$\varepsilon$.
It is straightforward then to check the required assumptions by using
\ref{LLambdasmall},
and then the results of
\cite[Appendix B]{kapouleas:annals} apply.
We upgrade the $L^2$ estimates for $f_0-1$
to $C^{2,\beta}$ estimates on $S_5[1]$
by using the uniformity of the geometry of $S_6[1]$ (see \ref{Lh})
and standard linear theory interior estimates.
Applying then the variant of \ref{CLambdau} we discussed earlier,
we estimate $f_0$ on $\Lambda$ and complete the proof.

\end{proof}

\subsection*{The (extended) substitute kernel}
$\phantom{ab}$
\nopagebreak

As we have already mentioned in the introduction,
the extended substitute kernel in this case is particularly
simple since it is one-dimensional.
This reflects the fact that the approximate kernel,
and hence the substitute kernel also,
are one-dimensional.
Moreover decay can be ensured by using the substitute kernel
and so no further extended substitute kernel is required.
Motivated by proposition
\ref{Pappker}
above we define
a function $w\in C^\infty_\sym(M)$
by requiring that on $M\cap\D$ it satisfies
\addtocounter{theorem}{1}
\begin{equation}
\label{Ew}
w:=\psi[m^{-1},2m^{-1}](r).
\end{equation}
For future reference we record the following:

\addtocounter{equation}{1}
\begin{lemma}
\label{Lw}
Given $E\in C^{0}_\sym(\Stilde[1] )$ there is a unique $\mu\in\R$
such that
$\frac{2\rho^2}{|A|^2+m^2} (E+\mu w)$
is $L^2(\Stilde[1],h)$-orthogonal to $f_0$,
where $f_0$ is the eigenfunction in 
\ref{Pappker}.
Moreover
$$
|\mu|\le C         \left\|\,   \frac{2\rho^2}{|A|^2+m^2} E:L^2_\sym(\Stilde[1],h)\,    \right\|.
$$
\end{lemma}

\begin{proof}
Using \ref{LLambdasmall} and \ref{Erho} we conclude that
$\frac1C\le\frac{2\rho^2}{|A|^2+m^2}\le C$
on the support of $w$ which together with \ref{Lh}
implies the result.
\end{proof}

To arrange the decay we define $v\in C^{\infty}_\sym(\Stilde[1])$
by
\addtocounter{theorem}{1}
\begin{equation}
\label{Ev}
v:=f_0+u,
\end{equation}
where $u$ is the solution to $\Lchi u=-\Lchi f_0+\mu' w$ on $\Stilde[1]$
with vanishing Dirichlet data on
$C[0]\subset\partial\Stilde[1]$,
where $\mu'\in\R$ is determined by the requirement (recall \ref{Lw})
that
$$
E':=
\frac{2\rho^2}{|A|^2+m^2} (-\Lchi f_0+\mu' w)=
\lambda_0 f_0 +\mu' \frac{2\rho^2}{|A|^2+m^2} w
$$
is $L^2(\Stilde[1],h)$-orthogonal to $f_0$.
Note that 
the equation on $\Stilde[1]$
is equivalent to
$\Lh u = E'$,
and hence the orthogonality condition together with \ref{Lw}
implies the existence of a unique $u$.
We record now the properties of $v$:

\addtocounter{equation}{1}
\begin{lemma}
\label{Lv}
$v$ satisfies the following:
\newline
(i).
$\Lchi v=\mu_v w $ on $\Stilde[1]$ for some $\mu_v\in\R$,
and therefore $\Lchi v=0$ on $\Lambda$.
\newline
(ii).
$v=0$ on 
$C[0]\subset\partial\Stilde[1]$.
\newline
(iii).
$\|v:C^{2,\beta}(\Stilde[1],\chi)\|
\le
C  $.
\newline
(iv).
$|\mu_v|\le C\varepsilon$.
\newline
(v).
$\|   v-1   :   C^{2,\beta}(C_1[1],\chi)\,   \|
\le C(b) \varepsilon$,
where $\varepsilon$ is as in \ref{Pappker}.
\end{lemma}

\begin{proof}
(i) and (ii) follow from the definitions.
Using \ref{Pappker} and \ref{Lw} we have that 
$$
\left\| \lambda_0 f_0 +\mu' \frac{2\rho^2}{|A|^2+m^2} w :L^2(\Stilde,h)  \right\|
\le
C\varepsilon,
$$
which together with interior $C^{2,\beta}$
estimates on $S_5[1]$ allows us to conclude
(iii),  (iv), and (v).
\end{proof}

\subsection*{Solving the linearized equation semi-locally}
$\phantom{ab}$
\nopagebreak

In this subsection we solve and estimate the linear equation on
the extended standard regions.
We can assume the inhomogeneous term $E$ to vanish on $\Lambda_1$,
because in the proof of \ref{PMsolution}
we use \ref{PLambdaEV}  
to solve for the part of the inhomogeneous term
which is supported there.
In the case of $\Stilde[1]$ we have nontrivial approximate kernel and
therefore we have to adjust the inhomogeneous term appropriately
by using $w$.
$w$ can also be used so that appropriate exponential decay
can be arranged for the solution:

\addtocounter{equation}{1}
\begin{lemma}
\label{LStilde1E}
There is a linear map
$$
\Rcal_{\Stilde[1]}:
\{E\in C^{0,\beta}_\sym(\Stilde[1]): E \text{ is supported on }S_{1}[1]\}
\to
C^{2,\beta}_\sym(\Stilde[1])\times\R,
$$
such that the following hold for
$E$ in the domain of $\Rcal_{\Stilde[1]}$ above
and $(\varphi,\mu)=\Rcal_{\Stilde[1]}(E)$:
\newline
(i).
$\Lchi \varphi=E\, + \mu  w$ on $\widetilde{S}[1]$.
\newline
(ii).
$\varphi$ vanishes on $C[0]\subset\partial \widetilde{S}[1]$
and satisfies appropriate Neumann boundary conditions on
$C_\partial\subset\partial \widetilde{S}[1]$
(recall \ref{ECpartial}).
\newline
(iii).
$|\mu| +
\Vert \varphi: C_\sym^{2,\beta}(\widetilde{S}[1],\chi) \Vert 
\le C(b,\beta)\, \Vert E: C_\sym^{0,\beta}({S}_1[1],\chi)\Vert .$
\newline
(iv).
$\Vert \varphi: C_\sym^{2,\beta}(\Lambda,\chi,e^{-\gamma\xover})\Vert 
\le C(b,\beta,\gamma)\, \Vert E: C_\sym^{0,\beta}({S}_1[1],\chi)\Vert .$
\newline
(v).
$\Rcal_{\Stilde[1]}$ depends continuously on $\zeta$.
\end{lemma}

\begin{proof}
We fix $b$ to be large enough so that \ref{Pappker}, \ref{Lw},
and \ref{CLambdauV} with $\varepsilon_1=1$ apply.
By applying \ref{Lw} and using that $E$ is supported on $S_1[1]$
we have $\mu_1$ such that $|\mu_1|
\le C(b)\, \Vert E: C_\sym^{0,\beta}({S}_1[1],\chi)\Vert $
and
$\frac{2\rho^2}{|A|^2+m^2} (E+\mu_1 w)$
is $L^2(\Stilde[1],h)$-orthogonal to $f_0$.
There is a unique solution $\varphi_1\in C^{2,\beta}_\sym(\Stilde[1])$ which
is $L^2(\Stilde[1],h)$-orthogonal to $f_0$,
vanishes on 
$C[0]\subset\partial \widetilde{S}[1]$,
and satisfies $\Lchi\varphi_1=E+\mu_1 w$ on $\Stilde[1]$.
Using interior estimates on $S_2[1]$ for $\varphi_1$ and
applying \ref{CLambdauV} on $\Lambda_{0,1}$ with $u=v-\avg v$
on $C_1[1]\subset\partial\Lambda_{0,1}$,
and once more with $u=\varphi_1-\avg\varphi_1$
on $C_1[1]\subset\partial\Lambda_{0,1}$,
we determine $\mu_2$
such that by taking $\varphi:=\varphi_1+\mu_2 v$,
and $\mu:=\mu_1+\mu_2\mu_v$
and using \ref{CLunique}
and the available estimates from 
\ref{CLambdauV} 
and \ref{Lv} we complete the proof.
\end{proof}

The corresponding statement for $\Stilde[0]$ is simpler,
reflecting the triviality
of the approximate kernel there and that we do not need exponential decay either:

\addtocounter{equation}{1}
\begin{lemma}
\label{LStilde0E}
There is a linear map
$$
\Rcal_{\Stilde[0]}:
\{E\in C^{0,\beta}_\sym(\Stilde[0]): E \text{ is supported on }S_{1}[0]\}
\to
C^{2,\beta}_\sym(\Stilde[0]),
$$
such that the following hold for
$E$ in the domain of $\Rcal_{\Stilde[0]}$ above
and $\varphi=\Rcal_{\Stilde[0]}(E)$:
\newline
(i).
$\Lchi \varphi=E$ on $\widetilde{S}[0]$.
\newline
(ii).
$\varphi$ vanishes on $\partial \widetilde{S}[0]$.
\newline
(iii).
$
\Vert \varphi: C_\sym^{2,\beta}(\widetilde{S}[0],\chi) \Vert 
\le C(b,\beta)\, \Vert E: C_\sym^{0,\beta}({S}_1[0],\chi)\Vert .$
\newline
(iv).
$\Rcal_{\Stilde[0]}$ depends continuously on $\zeta$.
\end{lemma}

\begin{proof}
By \ref{Pappker} there are no small eigenvalues and so we can solve
and obtain $L^2(h)$ estimates which together with interior estimates on $S_2[0]$
and \ref{CLambdau} imply the result.
\end{proof}

\subsection*{Solving the linearized equation globally}
$\phantom{ab}$
\nopagebreak

In order to solve the linearized equation \ref{Elinear}
globally on $M$ and provide estimates for the solutions,
we paste together the semi-local solutions provided by
\ref{PLambdaEV},
\ref{LStilde1E},
and
\ref{LStilde0E}
to obtain a global solution in the proof of \ref{PMsolution}.
Before we state the Proposition we define appropriate norms:

\addtocounter{equation}{1}
\begin{definition}
\label{Dnorm}
For $k\in\mathbb{N}$
and
$\beta,\gamma\in(0,1)$
we define a norm $\|\, .\,\|_{k,\beta,\gamma}$
on $C^{k,\beta}_\sym(M)$
by 
$$
\|\phi\|_{k,\beta,\gamma}:=\|\phi:C^{k,\beta}_\sym(M,\chi,\widetilde{f})\|,
$$
where the weight function $\widetilde{f}$ is defined by requesting that it is invariant
under the action of $\grouptilde$,
$\widetilde{f}=1$ on $S[1]$,
$\widetilde{f}=e^{-\gamma\xover}$ on $\Lambda$,
and 
$\widetilde{f}=e^{(\azero-2b)\gamma}=\left.e^{-\gamma\xover}\right|_{C[0]}$ on $S[0]$
(recall \ref{Eell}).
\end{definition}

Note that $\widetilde{f}$ is continuous and its minimum
as well the maximum of $\rho\widetilde{f}$ are attained
on $S[0]$, and therefore
using \ref{Ea} we have
\addtocounter{theorem}{1}
\begin{equation}
\label{Ef}
\tb^{\frac89\gamma+\frac19}\le\widetilde{f}
\quad\text{ and }\quad
\rho\widetilde{f}\le\tb^{\frac89\gamma-1}\qquad\text{on}\quad M.
\end{equation}

Before we proceed to state and prove the main Proposition of this section,
we give an estimate of the inhomogeneous term $E=\rho^{-2} H$ of the main linearized
equation in this paper:

\addtocounter{equation}{1}
\begin{lemma}
\label{LEH}
If $m$ is large enough in terms of $\gamma$ we have on $M$ the estimate
$$
\|\rho^{-2}H\|_{2,\beta,\gamma}\le
C\tau.
$$
\end{lemma}

\begin{proof}
Using 
\ref{Ea}, \ref{EMhattor}, \ref{Eacat}, \ref{Erho}, and \ref{Ef},
we easily check that
$|\zz|+\tau\le
m^2\,\tau$,
$m^2\,\tau^2\le
 \widetilde{f} \tau$,
and
$\rho^{-2}\,m^2\tau\le
C \widetilde{f} \tau$.
These imply that
$(\tau+\rho^{-2})(|\zz|+\tau)\le
C \widetilde{f} \tau$,
which by \ref{Lquantities} implies the result.
\end{proof}

\addtocounter{equation}{1}
\begin{prop}
\label{PMsolution}
There is a linear map $\Rcal_M:C^{0,\beta}_{sym}(M)\to C^{2,\beta}_{sym}(M)\times\R$
such that for $E \in C^{0,\beta}_{sym}(M)$
and $(\varphi,\mu)=\Rcal_M E$
the following hold:
\newline
(i).
$\Lchi \varphi=E\, + \mu  w$ on $M$.
\newline
(ii).
$|\mu|+\|\varphi\|_{2,\beta,\gamma}\le C(b,\beta,\gamma) \, \|E\|_{0,\beta,\gamma}$.
\newline
(iii).
$\Rcal_M$ depends continuously on $\zeta$.
\end{prop}

\begin{proof}
We decompose $E=E_{S[0]}+E_{S[1]}+E_\Lambda$
by requesting that 
$E_{S[0]}$, $E_{S[1]}$, and $E_\Lambda$,
are invariant under $\grouptilde$ and satisfy
$$
\begin{aligned}
E_{S[0]}:=&E\,\psi[1,0]\circ \xunder,\\
E_{S[1]}:=&E\,\psi[1,0]\circ \xover,\\
E_\Lambda :=&E\,\psi[0,1]\circ \xboth,
\end{aligned}
$$
on $\Lambda$,
$ E_{S[0]}:=E$,
$E_{S[1]}:=0$,
$E_\Lambda :=0$ on $S[0]$,
and 
$ E_{S[0]}:=0$,
$E_{S[1]}:=E$,
$E_\Lambda :=0$ on $S[1]$.
Using \ref{PLambdaEV} we define $V_\Lambda\in C^{2,\beta}_\sym(M)$
by
$V_\Lambda=0$ on $S[0]\cup S[1]$
and $V_\Lambda=\psi[0,1]\circ\xboth  \,
\Rcal_\Lambda\, E_\Lambda $ on $\Lambda$.
$\Lchi V_\Lambda-E_\Lambda$ is supported
on $\Lambda\setminus\Lambda_1$,
and can be decomposed as 
$\Lchi V_\Lambda-E_\Lambda=\underline{E}+\overline{E}$
where $\underline{E}$ is supported on $\{\xunder\le1\}$
and $\overline{E}$ is supported on $\{\xover\le1\}$.

Using \ref{LStilde1E} and \ref{LStilde0E}
we define 
$V_{S[0]}\in C^{2,\beta}_\sym(M)$
and
$V_{S[1]}\in C^{2,\beta}_\sym(M)$
by requesting the following:
$V_{S[1]}=0$ on $S[1]$
and
$V_{S[1]}=\psi[0,1]\circ\xunder\, V'_{S[1]}$
on $\Stilde[1]$,
where
$$
(V'_{S[1]},\mu_1)=\Rcal_{\Stilde[1]}(E_{S[1]}-\overline{E}).
$$
$V_{S[0]}=0$ on $S[1]$ and
$V_{S[0]}=
\psi[0,1]\circ\xover
\,
\Rcal_{\Stilde[0]}(E_{S[0]}-\underline{E})$
on $\Stilde[0]$.
We define then $\varphi_1:=V_\Lambda+V_{S[0]}+V_{S[1]}$ and $E_1$ by
$\Lchi\varphi_1+E_1=E+\mu_1 w$.
We iterate with $E_1$ instead of $E$ and so on.
We define then $\varphi:=\sum_{n=1}^\infty\varphi_n$
and $\mu:=\sum_{n=1}^\infty\mu_n$
and complete the proof by using the estimates and results of
\ref{PLambdaEV},
\ref{LStilde1E},
and \ref{LStilde0E},
where 
\ref{LStilde1E}
is applied with $\gamma'=\frac{\gamma+1}2$ in place of $\gamma$.
\end{proof}

\section{The main results}
\label{Sresults}
\nopagebreak

\subsection*{The nonlinear terms}
$\phantom{ab}$
\nopagebreak

If $\phi\in C^1_\sym(M)$ is appropriately small,
we denote by $M_\phi$ the perturbation of $M$ by $\phi$,
defined as $I_\phi(M)$ in the notation of Appendix \ref{AppendixB},
where $I:M\to\Sph^3(1)$ is the inclusion map of $M$.
Clearly then $M_\phi$ is invariant under the action of $\grouptilde$
on the sphere $\Sph^3(1)$.
Using then rescaling and Proposition \ref{PXpert}
we prove a global estimate of the nonlinear terms
for the mean curvature of $M_\phi$ as follows:

\addtocounter{equation}{1}
\begin{lemma}
\label{Lquadratic}
If $\phi\in C^{2,\beta}_\sym(M)$
satisfies
$\|\phi\|_{2,\beta,\gamma}<\tb^{1-\frac{3\gamma}4}$,
then $M_\phi$ is well defined as above and satisfies
$$
\|\rho^{-2} H_\phi-\rho^{-2}H-\Lchi\phi\|_{0,\beta,\gamma}
\le
\tb^{\frac{3\gamma}4-1}\|\phi\|_{2,\beta,\gamma}^2,
$$
where
$H_\phi$ is the mean curvature of $M_\phi$
(pulled back to $M$ by $I_\phi$),
and $H$ is the mean curvature of $M$.
\end{lemma}

\begin{proof}
Let $D$ be a disc of radius $1$ and center at some point $p\in M$ with respect to the $\chi$
metric.
If we magnify the metric of the sphere $\Sph^3(1)$ by a factor $\rho(p)$
it is easy to arrange for the hypothesis \ref{Ec1} to be satisfied
so that we can apply \ref{PXpert} with some universal $c_1$ to conclude
$$
\|(\rho(p))^{-1}(H_\phi-H-\Lcal\phi):C^{0,\beta}(D,\chi)\|
\le
\frac1{\epsilon(c_1)}\|\rho(p)\phi:C^{2,\beta}(D,\chi)\|^2,
$$
where the factors $\rho(p)$ correspond to the scaling of the quantities involved.
By the multiplicative properties of the Holder norms we conclude
$$
\|\rho^{-2}(H_\phi-H-\Lcal\phi):C^{0,\beta}(D,\chi)\|
\le
\frac{\rho(p)}{\epsilon(c_1)}\|\phi:C^{2,\beta}(D,\chi)\|^2.
$$
By \ref{Dnorm} we conclude
$$
\frac1{\widetilde{f}(p)}\|\rho^{-2}(H_\phi-H-\Lcal\phi):C^{0,\beta}(D,\chi)\|
\le
\frac{\rho(p)\widetilde{f}(p)}{\epsilon(c_1)}\|\phi\|^2_{2,\beta,\gamma}.
$$
This implies the result by using \ref{Ef}.
\end{proof}

\subsection*{The vertical force and balancing}
$\phantom{ab}$
\nopagebreak

If $\phi\in C^{1}_\sym(M)$, $M_\phi$, and $H_\phi$
are as in the previous subsection we define $\Fcal$ by
\addtocounter{theorem}{1}
\begin{equation}
\label{EFcal}
\Fcal:=\int_{M_\phi\cap\D_+} H_\phi\,\left<\nu,\vecK\right>\,dg
=\int_{M_\phi\cap\partial\D_+} \left<\veceta,\vecK\right>\,dg,
\end{equation}
where $\D_+:=\D\cap\{\zz\ge0\}$,
$\nu$ the unit normal chosen so that $\left<\nu,\partial_\zz\right>   >0$ on $\Mhattor$,
$\vecK$ is the Killing field defined in \ref{EvecK},
and $\eta$ the outward conormal to $\partial(M_\phi\cap\D_+)=M_\phi\cap\partial\D_+$
tangent to $M_\phi$.
Note that the second equality in 
\ref{EFcal}
follows from the first variation formula \cite{Si,KKS}.
We have then the following, where we could be using $\|\phi\|_{1,0,\gamma}$
instead of $\|\phi\|_{2,\beta,\gamma}$ as well:

\addtocounter{equation}{1}
\begin{lemma}
\label{LFcal}
If $\|\phi\|_{2,\beta,\gamma}<\tb^{1-\frac\gamma4}$, then there is a universal constant $C$ such that
$$
\left| \frac{m^2}{8\tau\pi^2}\Fcal+\zeta\right|
\le C\,(1+\frac1\tau\|\phi\|_{2,\beta,\gamma}).
$$
\end{lemma}

\begin{proof}
Let $d:=\piover$ and decompose 
$$
\partial(M_\phi\cap\D_+)=M_\phi\cap\partial\D_+=
\partial_{+1}\cup\partial_{-1}\cup\partial_{+2}\cup\partial_{-2}\cup\partial_0,
$$
where
$\partial_{+1}\subset\{\xx=d\}$,
$\partial_{-1}\subset\{\xx=-d\}$,
$\partial_{+2}\subset\{\yy=d\}$,
$\partial_{-2}\subset\{\yy=-d\}$,
and
$\partial_{0}\subset\{\zz=0\}$.
We use the big-$O$ notation to denote terms $O(A)$ which satisfy $|O(A)|\le C A$
for some universal constant $C$.
Using then \ref{EPhig} and \ref{EK} we calculate that on $\partial_{\pm1}$
\begin{equation*}
\begin{aligned}
\veceta&=\pm(1+\sin 2\zz)^{-1/2}\partial_\xx,
\\
\left<\veceta,\vecK\right>
&=-\frac1{\sqrt{2}}\sqrt{1+\sin2\zz}\,\,\cot(z+\frac\pi4)\,\,\sin\sqrt{2}\,d\,\,\cos \sqrt{2}\,\yy,
\\
dg&=\sqrt{1-\sin2\zz+\phi^2_\yy}\,d\yy.
\end{aligned}
\end{equation*}
Combining the above we obtain
$$
\int_{\partial_{\pm1}} \left<\veceta,\vecK\right>\,dg
=
-\frac1{\sqrt{2}}\int_{-d}^d(1-\sin2\zz +O(|z|^2+|\phi_\yy|^2))\,\,\sin\sqrt{2}\,d\,\,\cos \sqrt{2}\,\yy\,\,d\yy.
$$
Similarly
$$
\int_{\partial_{\pm2}} \left<\veceta,\vecK\right>\,dg
=
\frac1{\sqrt{2}}\int_{-d}^d(1+\sin2\zz +O(|z|^2+|\phi_\yy|^2))\,\,\cos \sqrt{2}\,\xx\,\,\sin\sqrt{2}\,d\,\,d\xx.
$$

Combining the above and using that on 
$\partial_{+1}\cup\partial_{-1}\cup\partial_{+2}\cup\partial_{-2}$
we have $\zz=\tau a+\phi$,
we conclude that
$$
\int_{\partial_{+1}\cup\partial_{-1}\cup\partial_{+2}\cup\partial_{-2}}
\left<\veceta,\vecK\right>\,dg
=(8a\tau+O(\|\phi\|+a^2\tau^2))  \,\,\, (2d^2+O(d^3)),
$$
where 
$\|\phi\|:=\|\phi\|_{2,\beta,\gamma}$.
Similarly
$$
\int_{\partial_{0}}
\left<\veceta,\vecK\right>\,dg
=
-2\pi\tau(1+O(\tau^{\frac\gamma2-1}\|\phi\|)).
$$
Combining and substituting $a=\frac{m^2}{4\pi}-\zeta+O(1)$ by \ref{Ea}
we conclude
$$
\Fcal=-\frac{8\pi^2\tau}{m^2}\zeta+\frac{1}{m^2}O(\tau+\|\phi\|)
,
$$
which implies the result.
\end{proof}

\subsection*{The main theorem}
$\phantom{ab}$
\nopagebreak

We have now all the information we need to state and prove the main theorem of the paper:

\addtocounter{equation}{1}
\begin{theorem}
\label{Tmain}
There are absolute constants $\cunder,C>0$ such that if $m$ is large enough,
then there is $\zeta_1\in [-\cunder,\cunder]$
such that on the corresponding initial surface $M$ there is $\phi\in C^\infty_\sym(M)$
with
$\|\phi\|_{2,\beta,\gamma}\le C\tau$
(with $\tau$ defined as in \ref{Etau})
such that $M_\phi$ is a genus $m^2+1$ embedded minimal surface
in $\Sph^3(1)$ invariant under the action of $\grouptilde$.
\end{theorem}

\begin{proof}
We will use a subscript $\zeta$ to specify the initial surface $M_\zeta$ which is
constructed as in the discussion preceding \ref{EM}.
We also define the map $\Xunder_\zeta:[-\azero,\azero]\times \Sph^1\to M_\zeta$
by requesting that $\Xunder_\zeta(\taa,\theta)=\Phi\circ X(a\taa/\azero,\theta)$
where the $X$ is the one defined for the given value of the parameter $\zeta$,
that is $\Xunder_\zeta$ is the parametrization corresponding to coordinates
$(\taa,\theta)$ for $\Mcat\subset M_\zeta$.
As in the proof of \ref{Lh} it is easy to check that there is
$\ttt:[\azero-2,\azero-1]\to[\azero-3,\azero]$
close to the identity map,
such that for $(\taa,\theta)\in[\azero-2,\azero-1]\times\Sph^1$
we have $\varpi\circ\Xunder_\zeta(\taa,\theta)=\varpi\circ\Xunder_0(\ttt(\taa),\theta)$.

We define now a diffeomorphism $F_\zeta:M_\zeta\to M_0$
by requiring that it is equivariant under the action of $\grouptilde$,
it satisfies $\varpi\circ F_\zeta=\varpi$ on $S_1[1]\subset M_\zeta$,
and that for
$(\taa,\theta)\in[-1,\azero]\times\Sph^1$
we have 
$$
F_\zeta\circ\Xunder_\zeta(\taa,\theta)=
\Xunder_0\left(\taa+\psi[\azero-2,\azero-1](\taa)\,(\ttt(\taa)-\taa)\,,\,\theta\right).
$$

We define now
a map $\Jcal : B \to B$ where
$$
B:=\{u\in C^{2,\beta}_\sym(M_0):\|u\|_{2,\beta,\gamma}\le\tb^{\frac\gamma2+1}\}
\times
[-\cunder,\cunder]
$$
as follows:
We assume $(u,\zeta)\in B$ given.
Let $\phi\in C^{2,\beta}(M_\zeta)$ be defined by
$\phi:=u\circ F_\zeta+\varphi$ where 
$(\varphi,\mu)=\Rcal_{M_\zeta} (-\rho^2 H)$
as in \ref{PMsolution}.
We have then
\newline
(a).
$\Lchi\varphi+\rho^{-2}H =\mu w_0$,
or equivalently
$\Lcal\varphi+H =\mu \rho^2 w_0$.
\newline
(b). By \ref{PMsolution} and \ref{LEH} we have
$$
|\mu|+\|\phi\|_{2,\beta,\gamma}\le C(b,\beta,\gamma) \, \tau.
$$

Applying 
\ref{PMsolution}
again and using \ref{Lquadratic}
we obtain
$(v,\mu'):=\Rcal_{M_\zeta} (-  (
\rho^{-2} H_\phi-\rho^{-2}H-\Lchi\phi) )$
which satisfies the following:
\newline
(c).
$\Lchi v+
\rho^{-2} H_\phi-\rho^{-2}H-\Lchi\phi 
=\mu' w_0$.
\newline
(d).
$|\mu'|+\|v\|_{2,\beta,\gamma}
\le
\tb^{\frac{3\gamma}4-1}\|\phi\|_{2,\beta,\gamma}^2$.

Combining (a) and (c) with the definition of $\phi$ we obtain
\newline
(e).
$\Lchi(v-u\circ F_\zeta)+\rho^{-2}H_\phi=(\mu+\mu')w_0$.

This motivates us to define
$$
\Jcal(u,\zeta)=\left(v\circ (F_\zeta)^{-1},
\frac{m^2}{8\tau\pi^2}\Fcal+\zeta\right),
$$
where $\Fcal$ is defined as in \ref{EFcal}.
By using (b), (d), and \ref{LFcal},
and by choosing $\cunder$ large enough in terms of an absolute constant,
it is straightforward to check that $\Jcal(B)\subset B$.
$B$ is clearly a compact convex subset of $C^{2,\beta'}_\sym(M_0)\times\R$
for $\beta'\in(0,\beta)$,
and it is easy to check that $\Jcal$ 
is a continuous map in the induced topology.
By Schauder's fixed point theorem
\cite[Theorem 11.1]{gilbarg} then,
there is a fixed point of $\Jcal$.
Using (e) then we conclude that for the corresponding $\zeta$ and $\phi$
we have 
$$
H_\phi=(\mu+\mu')\rho^2w_0,\qquad\qquad\Fcal=0.
$$
Since $\left<\nu,\vecK\right>  >0$ on the support of $w_0$ in
$(M_\zeta)_\phi$
the second equation implies that $\mu+\mu'=0$ and hence 
$(M_\zeta)_\phi$
is a minimal surface.
The smoothness of $\phi$
follows then by standard regularity theory.
The embeddedness of 
$(M_\zeta)_\phi$
follows from the smallness of $\|\varphi\|_{2,\beta,\gamma}$
and the size (by \ref{Ea}) of $a\tau$.
\end{proof}

\appendix

\section{A coordinate system on $\mathbb S^3(1)$}
\label{AppendixA}
\nopagebreak

\subsection*{The parametrization $\Phi$}
$\phantom{ab}$
\nopagebreak

It is very helpful that there is a coordinate system which is ideally suited
to describing the Clifford torus and its parallel surfaces.
We proceed to describe this coordinate system and the local parametrization
which is its inverse.
To simplify the notation we identify $\R^4\simeq\C^2\supset\Sph^3(1)$.
We define the parametrization $\Phi$,
which covers the unit sphere with two orthogonal circles removed,
that is
$\Sph^3(1)\setminus\{(z_1,z_2)\in \C^2:z_1=0\text{ or }z_2=0\}$,
by the following:
\addtocounter{theorem}{1}
\begin{equation}
\label{EPhi}
\begin{gathered}
  \Phi :  \domPhi
  \rightarrow   \mathbb S^3(1) \subset \mathbb \R^4\sim\C^2, 
\quad\text{ where}\quad
\domPhi:=
\mathbb R \times \mathbb R \times \left(-\frac{\pi}{4},\frac{\pi}{4}\right), \\
  \Phi(\xx,\yy,\zz) = 
    \cos(\zz + \tfrac{\pi}{4})
     \, e^{\sqrt{2\,}\yy  i\,}
      \vec{e}_1 
                        \,+\,
    \sin(\zz + \tfrac{\pi}{4})
     \, e^{\sqrt{2\,}\xx i\,}
      \vec{e}_2,
\end{gathered}
\end{equation}
where $\vec{e}_1=(1,0)$ and $\vec{e}_2=(0,1)$
form the standard basis of $\C^2$.

\subsection*{Symmetries of $\Phi$}
$\phantom{ab}$
\nopagebreak

To study the symmetries of the parametrization $\Phi$,
we first define for $c\in\R$ translations
$\XXX_c$, $\YYY_c$,
and reflections $\xbar_c$, $\ybar_c$, $\xbar:=\xbar_0$, $\ybar:=\ybar_0$, and $\zbar$,
of its domain $\domPhi$, by 
\addtocounter{theorem}{1}
\begin{equation}
\label{Edomainsymmetries}
\begin{aligned}
\XXX_c(\xx,\yy,\zz)&=(\xx+c,\yy,\zz),\qquad
&
\YYY_c(\xx,\yy,\zz)&=(\xx,\yy+c,\zz),
\\
\xbar_c(\xx,\yy,\zz)&=(2c-\xx,\yy,\zz),\qquad
&
\ybar_c(\xx,\yy,\zz)&=(\xx,2c-\yy,\zz),
\\
\zbar(\xx,\yy,\zz)&=(\yy,\xx,-\zz).
&
&
\end{aligned}
\end{equation}

We also define corresponding
rotations $\xtilde_c$, $\ytilde_c$,
and reflections
$\xbartilde_c$, $\ybartilde_c$, $\xbartilde:=\xbartilde_0$, $\ybartilde:=\ybartilde_0$,
and $\zbartilde$
of $\Sph^3(1)\subset\C^2$ by
\addtocounter{theorem}{1}
\begin{equation}
\label{Espheresymmetries}
\begin{aligned}
\xtilde_c(z_1,z_2)&=(z_1,
     \, e^{\sqrt{2}\,c\,i\,}
z_2),\qquad
&
\ytilde_c(z_1,z_2)&=(
     \, e^{\sqrt{2}\,c\,i\,}
z_1,z_2),\qquad
\\
\xbartilde(z_1,z_2)&=(z_1,\overline{z_2}),\qquad
&
\ybartilde(z_1,z_2)&=(\overline{z_1},z_2),\qquad
\\
\xbartilde_c&:=\xtilde_{2c}\circ\xbartilde,
&
\ybartilde_c&:=\ytilde_{2c}\circ\ybartilde,
\\
\zbartilde(z_1,z_2)&=(z_2,z_1).\qquad
&&
\end{aligned}
\end{equation}

Note that $\xbartilde_c$, $\ybartilde_c$ and $\zbartilde$
are reflections with respect to the $3$-planes
$\left< \vec{e}_1, i \vec{e}_1 ,  e^{\sqrt{2}\,c\,i\,} \vec{e}_2 \right>_\R$,
$\left< e^{\sqrt{2}\,c\,i\,} \vec{e}_1, \vec{e}_2, i \vec{e}_2 \right>_\R$,
and the 2-plane $\{z_1=z_2\}$ respectively.
$\zbartilde$ exchanges the two sides of the Clifford torus and also
interchanges its parallels with its meridians.
$\xtilde_{\sqrt{2\,}\,\pi}$
and
$\ytilde_{\sqrt{2\,}\,\pi}$
are the identity map.
We record the symmetries of $\Phi$ in the following lemma:

\addtocounter{equation}{1}
\begin{lemma}
\label{LPhisymmetries}
$\Phi$ is a covering map onto 
$\Sph^3(1)\setminus\{(z_1,z_2)\in \C^2:z_1=0\text{ or }z_2=0\}$.
Moreover the following hold:
\newline
(i).
The group of covering transformations is generated by
$\XXX_{\sqrt{2\,}\,\pi}$ and $\YYY_{\sqrt{2\,}\,\pi}$,
in particular
$\Phi=\Phi\circ\XXX_{\sqrt{2\,}\,\pi} =\Phi\circ\YYY_{\sqrt{2\,}\,\pi}$.
\newline
(ii).
$\xbartilde_c\circ\Phi=\Phi\circ\xbar_c$,
$\ybartilde_c\circ\Phi=\Phi\circ\ybar_c$,
and
$\zbartilde\circ\Phi=\Phi\circ\zbar$.
\newline
(iii).
$\xtilde_c\circ\Phi=\Phi\circ\XXX_c$
and
$\ytilde_c\circ\Phi=\Phi\circ\YYY_c$.
\end{lemma}

\begin{proof}
(ii) and (iii) follow from the definitions.
(i) follows from (iii) and the observation that 
$\xtilde_{\sqrt{2\,}\,\pi}$
and
$\ytilde_{\sqrt{2\,}\,\pi}$
are the identity map.
\end{proof}

\subsection*{The coordinates $\xx\yy\zz$}
$\phantom{ab}$
\nopagebreak

The local inverses of $\Phi$ provide us with local coordinate systems.
We denote the corresponding coordinates by $\xx,\yy,\zz$.
A straightforward calculation shows that
\addtocounter{theorem}{1}
\begin{equation}
\label{Epartialxyz}
\begin{aligned}
\partial_\xx&=
    \sqrt{2\,}\,\sin(\zz + \tfrac{\pi}{4})
     \, i \, e^{\sqrt{2}\,\xx\,i\,}
      \vec{e}_2,\\
\partial_\yy&=
    \sqrt{2\,}\,\cos(\zz + \tfrac{\pi}{4})
     \, i \, e^{\sqrt{2}\,\yy\,i\,}
      \vec{e}_1 ,\\
\partial_\zz&=
   - \sin(\zz + \tfrac{\pi}{4})
     \, e^{\sqrt{2}\,\yy\,i\,}
      \vec{e}_1 
                        \,+\,
    \cos(\zz + \tfrac{\pi}{4})
     \, e^{\sqrt{2}\,\xx\,i\,}
      \vec{e}_2.
\end{aligned}
\end{equation}
By calculating further we obtain
\addtocounter{theorem}{1}
\begin{equation}
\label{EPhig}
    \Phi^* g = (1+\sin2\zz)\,d\xx^2 + (1-\sin2\zz)\,d\yy^2 + d\zz^2,
\end{equation}
where $g$ is the induced metric on the unit sphere $\Sph^3(1)$.
Moreover the only non-vanishing Christoffel symbols for the $(\xx,\yy,\zz)$-coordinate
system are given by
\addtocounter{theorem}{1}
\begin{equation}
\label{EChr}
\begin{aligned}
\Gamma^1_{13}&=\frac{\cos2\zz}{1+\sin2\zz},
\qquad&
\Gamma^2_{23}&=-\frac{\cos2\zz}{1-\sin2\zz},
\\
\Gamma^3_{11}&=-\cos2\zz,
\qquad&
\Gamma^3_{22}&=\cos2\zz.
\end{aligned}
\end{equation}

The level surface with $\zz=0$ is
the Clifford torus
\addtocounter{theorem}{1}
\begin{equation}
\label{ECT}
\T:=\Phi(\{\zz=0\})=\{(z_1,z_2)\in\Sph^3(1)\subset\C^2:|z_1|=|z_2|=1/\sqrt{2\,}\}.
\end{equation}
The level surfaces $\Phi(\{\zz=c\})$ ($\zz\in(-\frac\pi4,\frac\pi4)$)
are tori of constant mean curvaure,
parallel at distance $c$ to the Clifford torus $\T$,
with $\partial_{\zz}$ as their unit normal vector field.
Note also that for $c\in\R$
we have the level surfaces 
\begin{equation*}
\begin{aligned}
\Phi(\{\xx=c\})=\{t_1\vece_1+t_2i\vece_1+t_3
     \, e^{\sqrt{2}\,c\,i\,}\vece_2:
t_1,t_2\in\R,t_3\in\R^+\}
\cap\,\Sph^3(1), \\
\Phi(\{\yy=c\})=\{
 t_1    \, e^{\sqrt{2}\,c\,i\,}\vece_1+t_2\vece_2+t_3i\vece_2
:
t_1\in\R^+,t_2,t_3\in\R\}
\cap\,\Sph^3(1), 
\end{aligned}
\end{equation*}
which 
are equatorial half-two-spheres orthogonal to the parallel tori.
These three families of level surfaces are orthogonal.
The intersections of the last two
are great semicircles orthogonal to the tori.
Finally a calculation shows that $ det[ \Phi,
\Phi_\xx, \Phi_\yy, \Phi_\zz ] = \cos2\zz >0$.

\subsection*{Killing fields}
$\phantom{ab}$
\nopagebreak

Clearly 
$ \partial_{\xx}$ and $ \partial_{\yy}$ are Killing fields generating the
rotations in the $\left< \vec{e}_2, i\vec{e}_2\right>_\R$  and  
$\left< \vec{e}_1, i\vec{e}_1\right>_\R$
planes respectively.
However
$\partial_{\zz}$ is not a Killing field.
For this reason we consider the Killing field $\vecK$ 
which agrees with $\partial_z$ at $\Phi(0,0,0)=2^{-1/2}(\vece_1+\vece_2)$
and is defined by
\addtocounter{theorem}{1}
\begin{equation}
\label{EvecK}
\left. \vecK\right|_{(z_1,z_2)}
:=
-\real z_2 \, \vece_1
+
\real z_1 \, \vece_2.
\end{equation}
$\vecK$ generates the rotations in the $\left<\vece_1,\vece_2\right>_\R$ plane.
A straightforward calculation shows that
\addtocounter{theorem}{1}
\begin{multline}
\label{EK}
\vecK=
    -\tfrac{1}{\sqrt{2}}\cot(\zz+
\pi/4) \sin\sqrt{2}\xx \cos\sqrt{2}\yy
    \partial_\xx  \\
\quad + \tfrac{1}{\sqrt{2}}\tan(\zz + \pi/4) \cos\sqrt{2}\xx
    \sin\sqrt{2}\yy \partial_\yy + \cos\sqrt{2}\xx \cos\sqrt{2}\yy \partial_\zz.
\end{multline}

\section{The mean curvature of a perturbed surface}
\label{AppendixB}
\nopagebreak

We assume given an immersion $X:D\to U$,
where $D$ is a disc of radius $1$ in the Euclidean plane $\R^2$,
and $U$ is an open cube in $\R^3$
equipped with a metric $g$ whose components are functions $g_{ij}:U\to \R$.
We assume that the following holds for some $c_1>0$:
\begin{equation}
\addtocounter{theorem}{1}
\label{Ec1}
\Vert \partial X:C^{2,\beta}(D,g_0)\Vert  \le c_1,
\qquad
\|g_{ij}:C^{2,\beta}(U,g_0)\| \le c_1,
\qquad
g_0 \le c_1 X^*g,
\end{equation}
where $\partial X$ are the partial derivatives of the coordinates of $X$, 
and $g_0$ denotes the standard Euclidean metric on $U$ or $D$ respectively.
Note that \ref{Ec1} can be arranged by 
first appropriately magnifying the target (see for example \ref{Lquadratic}).
We also choose a unit normal $\nu:D\to \R^3$ for the immersion $X$
with respect to the $g$ metric.
Given a function $\phi:D\to \R$ which is small enough
we define $X_\phi:D\to U$ by
\addtocounter{theorem}{1}
\begin{equation}
\label{EXphi}
X_\phi(p):=\exp_{X(p)}(\phi(p)\,\nu(p)),
\end{equation}
where $\exp$ is the exponential map with respect to the $g$ metric.
We have then the following:

\begin{prop}
\addtocounter{equation}{1}
\label{PXpert}
There exists a (small) constant $\epsilon(c_1)>0$
such that if $X$ is an immersion satisfying \ref{Ec1}
and the function $\phi:D \ra \R$ satisfies
$$
\Vert \phi:C^{2,\beta}(D,g_0)\Vert  < \epsilon(c_1),
$$
then $X_\phi : D \ra U$ is a well-defined immersion by \ref{EXphi}
and satisfies
$$
\Vert
H_\phi-H-(\Delta_g+|A|^2+Ric(\nu,\nu))\phi:
C^{0,\beta}(D,g_0)\|
\le
\frac1{\epsilon(c_1)}
 \Vert \phi:C^{2,\beta}(D,g_0)\Vert^2,
$$
where $H=\tr_gA$ is the mean curvature of $X$, defined as the trace with
respect to $X^*g$ of the second fundamental form $A$,
$H_\phi$ is the mean curvature of $X_\phi$,
$\Delta_g$ is the Laplacian with respect to $X^*g$,
and $Ric$ is the Ricci curvature of $(U,g)$.
\end{prop}

\begin{proof}
That the linear terms are as stated is well known and follows by a straightforward calculation we omit.
The nonlinear terms are given by expressions of monomials consisting of contractions
of derivatives of $X$ and derivatives of
$\phi$.
This implies both the existence results and the estimate on the nonlinearity.
\end{proof}

\bibliographystyle{amsplain}
\bibliography{paper}
\end{document}